\documentclass[a4paper, 11pt, twoside, openany]{article} 
\usepackage{amsmath,amssymb,stmaryrd, mathtools}
\usepackage[all]{xy}
\usepackage{fancyhdr}
\usepackage{footmisc}
\usepackage{hyperref}
\usepackage[margin=22mm]{geometry}
\usepackage[UKenglish]{datetime}
\usepackage{graphicx}
\usepackage{wrapfig}
\usepackage[usenames, dvipsnames]{color}
\usepackage{accents}
\usepackage{cases}
\usepackage{enumerate}
\usepackage{tikz-cd}

\usepackage[normalem]{ulem}

\newcommand{\bs}[1]{\boldsymbol{#1}}

\newcommand{\pp}[2]{\frac{\partial #1}{\partial #2}} 
\newcommand{\dd}[2]{\frac{\delta #1}{\delta #2}}

\numberwithin{equation}{section}

\newcounter{savefootnote}
\newcounter{symfootnote}
\newcommand{\symfootnote}[1]{%
   \setcounter{savefootnote}{\value{footnote}}%
   \setcounter{footnote}{\value{symfootnote}}%
   \ifnum\value{footnote}>8\setcounter{footnote}{0}\fi%
   \let\oldthefootnote=\thefootnote%
   \renewcommand{\thefootnote}{\fnsymbol{footnote}}%
   \footnote{#1}%
   \let\thefootnote=\oldthefootnote%
   \setcounter{symfootnote}{\value{footnote}}%
   \setcounter{footnote}{\value{savefootnote}}%
}

\begin{document}
\nocite{*}
%
%
\begin{center}
{\Large Structure preserving transport stabilized compatible finite element methods for magnetohydrodynamics}
\end{center}
\vspace{-3mm}
\hrulefill
\begin{center}
{Golo A. Wimmer$^{1}$\symfootnote{correspondence to: gwimmer@lanl.gov}, Xianzhu Tang$^1$}\\
\vspace{2mm}
{\textit{$^1$Los Alamos National Laboratory}}\\
\vspace{4mm}
\today
\end{center}
\vspace{5mm}
\begin{center}
\textbf{Abstract}
\end{center}
We present compatible finite element space discretizations for the ideal compressible magnetohydrodynamic equations. The magnetic field is considered both in div- and curl-conforming spaces, leading to a strongly or weakly preserved zero-divergence condition, respectively. The equations are discretized in space such that transfers between the kinetic, internal, and magnetic  energies are consistent, leading to a preserved total energy. We also discuss further adjustments to the discretization required to additionally achieve magnetic helicity preservation. Finally, we describe new transport stabilization methods for the magnetic field equation which maintain the zero-divergence and energy conservation properties, including one method which also preserves magnetic helicity. The methods' preservation and improved stability properties are confirmed numerically using a steady state and a magnetic dynamo test case.
\section{Introduction}
The compressible magnetohydrodynamic (MHD) equations are of interest in a diverse range of application areas, including astrophysics and magnetic confinement fusion \cite{goedbloed2004principles}. Their most fundamental variant is given by single fluid ideal MHD, where dissipative effects are ignored and all terms except of the transport one are omitted in Ohm's law, which defines the electric field to be used in the magnetic field equation \cite{goedbloed2004principles, jardin2010computational}. While many physically relevant effects -- such as those related to Hall MHD -- are lost in the ideal case, the resulting system of partial differential equations still contains a rich structure that is challenging to discretize \cite{goedbloed2010advanced, jardin2010computational}. This structure in particular includes the magnetic field's zero-divergence property, and a failure to satisfy this property after discretization generally leads to large amounts of spurious noise. A wealth of numerical methods has been developed to avoid this problem, by means of, for example, Lagrange multipliers~\cite{schotzau2004mixed}, a post-processing step~\cite{brackbill1980effect}, 8-wave methods~\cite{powell1994approximate}, or discretizations that satisfy the zero-divergence property exactly~\cite{evans1988simulation}, \cite{hu2017stable} (and references therein). Another structure is given by a balance between the kinetic, internal, and magnetic energies, leading to conservation of total energy. In the context of hydrodynamics, it has long been known that a discrete analogue of this property may lead to more accurate long term predictions \cite{arakawa1997computational}. Furthermore, such a discrete analogue may lead to more accurate results at lower resolutions \cite{wimmer2020SUPG}. Additionally, the ideal MHD equations give rise to conserved helicities, including the magnetic helicity \cite{hu2021helicity} (and references therein). Finally, next to these quantities arising from the system's structure, the MHD equations also contain transport terms for each of the prognostic variables. When physical dissipative effects are relatively small -- which is often the case for many applications of interest -- additional stabilization methods may be needed to avoid spurious oscillations due to the discrete transport terms.\\ \\
%
Recently, a number of compatible finite element based discretizations have been described for various forms of the MHD equations, ensuring some or all of the aforementioned structure preserving properties \cite{gawlik2022finite, hu2021helicity, hu2017stable, laakmann2022structure}. In particular, the magnetic field's zero-divergence property is ensured to hold true via discrete vector calculus identities that are satisfied in the compatible finite element approach. Further, energy conservation is achieved by using matching discretizations for the various terms appearing in the system of equations. This approach may be facilitated by a rich underlying theory starting from a Lagrangian mechanics point of view \cite{gawlik2022finite}, or a Hamiltonian one based on Poisson brackets \cite{morrison1982poisson}. Lastly, helicity conservation can be achieved by the use of appropriate projections, moving between the corresponding div- and curl-conforming finite element spaces \cite{hu2021helicity}. In the context of compatible finite element discretizations for hydrodynamics in ocean and atmosphere modeling, transport stabilization methods have further been incorporated without compromising on energy conservation for the density, thermal, velocity, and vorticity fields \cite{bauer2018energy, eldred2019quasi, lee2021petrov, natale2018variational, wimmer2020SUPG, wimmer2020density, wimmer2020energy}. However, there is a lack of such methods specifically aimed at improving the magnetic field transport term's stability, which comes with the added difficulty of maintaining the latter field's zero-divergence property. One recent approach to resolve this is given by \cite{wu2020simplex}, which concerns general grad-, div-, and curl-conforming convection-diffusion problems, and relies on exponential fitting methods.\\ \\
In this paper, we introduce two new transport stabilization methods for the div- and curl-conforming discretizations of the magnetic field equation. They are based on interior penalty type formulations \cite{burman2005unified, burman2006edge}, and can be seen as sub-grid resistivity terms. Unlike the addition of a small physical resistive contribution for stabilization purposes, the latter two methods are consistent with the ideal MHD equations in the sense that for a continuous current density, the interior penalty terms vanish. The two methods differ in that the first one's penalty term is analogous to the resistivity term, while in the second one, the penalty is adjusted to vanish in the computation of the rate of change of total helicity. Further, we examine a residual based SUPG formulation for the magnetic field equation, which in the context of compatible finite element methods has been considered for atmosphere and ocean modeling for the vorticity equation before \cite{bauer2018energy, wimmer2020energy}. All three methods preserve the zero-divergence as well as energy conservation properties. Finally, the three methods' stabilization properties are studied numerically and compared against each other for both a curl- and a div-conforming discretization of the magnetic field evolution equation.\\ \\
While magnetic field transport constitutes the main focus of the paper, we also introduce new energy conserving transport stabilization methods for the thermal field. Transport stabilized thermohydrodynamic systems have been studied using the compatible finite element framework before (e.g. \cite{eldred2019quasi} for buoyancy and \cite{wimmer2020SUPG} for potential temperature); however, we found alternative approaches to fit better into the context of the compressible MHD equations considered here. In particular, we include an energetically neutral interior penalty term when using the temperature as the thermal field, and an energetically neutral upwind formulation if pressure is used instead.  \\ \\
The remainder of this work is structured as follows: In Section \ref{sec_background}, we review the compressible magnetohydrodynamic equations and their structure preserving properties, as well as the compatible finite element method. In Section \ref{sec_Space_discretization}, we introduce the structure preserving space discretization, including stabilization methods for all transport terms. In Section \ref{sec_Numerical_results}, we present and discuss numerical results. Finally, in Section \ref{sec_conclusion}, we review our results and discuss ongoing work.
%
\section{Background} \label{sec_background}
In this section, we introduce the compressible magnetohydrodynamic equations, and review some of its structure preserving properties. Further, we briefly review the compatible finite element spaces to be used in the space discretization.
%
\subsection{Compressible MHD}
In this work, we derive discretizations for the ideal compressible magnetohydrodynamic equations. However, in order to motivate the new interior penalty based stabilization methods to be introduced in the next section, we will start from the resistive MHD equations. They are given by
\begingroup
\addtolength{\jot}{2mm}
\begin{subequations} \label{MHD_cts}
\begin{align}
&\pp{n}{t} + \nabla \cdot (n \mathbf{V}) = 0, \label{cty_eqn_full}\\
&\pp{\mathbf{V}}{t} + (\nabla \times \mathbf{V}) \times \mathbf{V} + \frac{1}{2} \nabla |\mathbf{V}|^2 + \frac{1}{m_i n }\nabla ( n T) - \frac{1}{m_i n }\mathbf{j} \times \mathbf{B} = 0, \label{momentum_eqn} \\
\frac{n}{\gamma - 1} &\left( \pp{T}{t} + \mathbf{V}\cdot \nabla T \right) +  n T \nabla \cdot \mathbf{V} = \eta |\mathbf{j}|^2 \label{T_eqn_cts}\\
&\pp{\mathbf{B}}{t} + \nabla \times \mathbf{E} = 0, \label{B_eqn_cts}\\
&\mathbf{E} = - \mathbf{V} \times \mathbf{B} + \eta \mathbf{j}, \hspace{1cm} \mathbf{j} = \frac{1}{\mu_0} \nabla \times \mathbf{B}, \label{Ohms_law}
\end{align}
\end{subequations}
\endgroup in a domain $\Omega$, and where the prognostic fields are
given by the ion number density $n$, ion flow velocity $\mathbf{V}$,
plasma temperature $T$, and magnetic field $\mathbf{B}$. Details of
the field's underlying spaces will be postponed to the next section,
and in this section, we simply assume sufficient continuity for the
differential operations occurring in the above equations to be
well-defined. Note that in view of this section's discussion on energy
conservation, we wrote the velocity transport term in its vector
invariant form,
\begin{equation}
(\mathbf{V} \cdot \nabla) \mathbf{V} = (\nabla \times \mathbf{V}) \times \mathbf{V} + \frac{1}{2} \nabla |\mathbf{V}|^2.
\end{equation}
Additionally, there are two diagnostic fields, given by the electric
field $\mathbf{E}$, and the current density $\mathbf{j}$. The
constants $m_i$ and $\mu_0$ denote the ion mass and vacuum
permeability, respectively, and $\gamma = 5/3$. Finally, $\eta$
denotes the plasma resistivity. If $\eta \equiv 0$ -- i.e. if
there is no resistive or Ohmic heating effect -- then the system of
equations reduces to ideal MHD.\\ \\
\textbf{Remark 1.} Alternatively to the temperature $T$, it is also possible to consider other thermal fields such as the pressure $p = nT$, with a governing equation given by.
\begin{equation}
\frac{1}{\gamma - 1} \left( \pp{p}{t} + \nabla \cdot (p \mathbf{V}) \right) +  p \nabla \cdot \mathbf{V} = \eta |\mathbf{j}|^2. \label{p_eqn_cts}
\end{equation}
This formulation leads to a more straightforward description of the system's energy transfers. However, here we consider the temperature $T$ instead, since we found our resulting discretization \eqref{nVT_discr} to be on average more accurate when compared to our pressure field based one \eqref{nVp_discr}. Further, the corresponding evolution equation for $T$ and its discretization are more amenable to incorporate an anisotropic heat flux -- which includes a gradient in $T$ rather than $p$ -- and which plays an important role in, for example, magnetic confinement fusion applications.\begin{center}\vspace{-15mm}\end{center}\hrulefill\\

Finally, we equip the above equations with free-slip boundary conditions for the velocity field, and mixed boundary conditions for the magnetic field. For a domain $\Omega$ with a boundary region $\partial \Omega$ and sub-regions $\partial \Omega_1$, $\partial \Omega_2$ such that $\partial \Omega_1 \cup \partial \Omega_2 = \partial \Omega$ and $\partial \Omega_1 \cap \partial \Omega_2 = \emptyset$, these are given by
\begin{align}
&\mathbf{V} \cdot \mathbf{n}|_{\partial \Omega} = 0, \\
&\mathbf{B} \cdot \mathbf{n}|_{\partial \Omega_1} = g_1, &\mathbf{H} \times \mathbf{n}|_{\partial \Omega_2} = \mathbf{g}_2 \label{B_bc},
\end{align}
for outward normal unit vector $\mathbf{n}$, and where $\mathbf{H} \coloneqq \mu_0 \mathbf{B}$. The magnetic field boundary conditions can alternatively also be expressed in terms of $\mathbf{E} \times \mathbf{n} = \mathbf{h}_1$ and $\mathbf{j} \cdot \mathbf{n} = h_2$ in $\partial \Omega_1$ and $\partial \Omega_2$, respectively \cite{bochev2002matching}. We note that this may be preferable depending on the choice of finite element function space, as will be shown in Section \ref{sec_Space_discretization}.
%
%
\\ \\
The above system of equations gives rise to a number of conserved quantities, and in the following, we will focus on the system's total energy, total magnetic helicity, and zero-divergence property for the magnetic field. The system's total energy is given by
\begin{equation}
H(n, \mathbf{V}, T, \mathbf{B}) =  \int_\Omega \left( \frac{1}{2} m_i n |\mathbf{V}|^2 + \frac{n T}{\gamma - 1} + \frac{|\mathbf{B}|^2}{2 \mu_0}\right) dx. \label{Hamiltonian}
\end{equation}
\textbf{Proposition 1a}. \textit{For boundary conditions of the form $\mathbf{V} \cdot \mathbf{n}|_{\partial \Omega} = 0$ and $\mathbf{B} \cdot \mathbf{n}|_{\partial \Omega} = 0$, the total energy is conserved in time.}\\ \vspace{-3mm} \\
\textbf{Proof.} The rate of change in time of the total energy \eqref{Hamiltonian} can be computed using the chain rule
\begin{equation}
\frac{dH}{dt} =  \left\langle \dd{H}{n}, \pp{n}{t} \right\rangle + \left\langle \dd{H}{\mathbf{V}}, \pp{\mathbf{V}}{t} \right\rangle + \left\langle \dd{H}{T}, \pp{T}{t} \right\rangle + \left\langle \dd{H}{\mathbf{B}}, \pp{\mathbf{B}}{t} \right\rangle, \label{Hamiltonian_full}
\end{equation}
where -- after some abuse of notation\footnote{Strictly speaking, the Hamiltonian derivatives are elements of the dual spaces of the respective fields, and $\langle ., . \rangle$ denotes a pairing of elements from a space and the space's dual.} -- $\langle ., . \rangle$ denotes the $L^2$ inner product, and the variational derivative related expressions are given by
\begin{align}
\dd{H}{n} = \frac{1}{2} m_i |\mathbf{V}|^2 +  \frac{T}{\gamma - 1}, \hspace{1cm} \dd{H}{T} = \frac{n}{\gamma - 1}, \hspace{1cm} \dd{H}{\mathbf{V}} = m_i n \mathbf{V}, \hspace{1cm} \dd{H}{\mathbf{B}} = \frac{\mathbf{B}}{\mu_0}. \label{H_variations}
\end{align}

Substituting these expressions and the prognostic fields' time derivatives in \eqref{Hamiltonian_full} by the evolution equations \eqref{cty_eqn_full} - \eqref{B_eqn_cts}, we obtain
\begingroup
\addtolength{\jot}{2mm}
\begin{subequations}
\begin{align}
&\left\langle \dd{H}{n}, \pp{n}{t} \right\rangle = - \left\langle \frac{1}{2} m_i |\mathbf{V}|^2, \nabla \cdot (n \mathbf{V}) \right\rangle - \left\langle \frac{T}{\gamma - 1}, \nabla \cdot (n \mathbf{V}) \right\rangle,\\
&\left\langle \dd{H}{\mathbf{V}}, \pp{\mathbf{V}}{t} \right\rangle = -\left\langle m_i n \mathbf{V}, (\nabla \times \mathbf{V}) \times \mathbf{V} \right\rangle - \left\langle m_i n \mathbf{V},  \frac{1}{2} \nabla |\mathbf{V}|^2 \right\rangle - \left\langle \mathbf{V}, \nabla ( n T)\right\rangle + \left\langle \mathbf{V}, \mathbf{j} \times \mathbf{B} \right\rangle, \label{dv_dvdt}\\
&\left\langle \dd{H}{T}, \pp{T}{t} \right\rangle = - \left\langle \frac{n}{\gamma - 1}, \mathbf{V} \cdot \nabla T\right\rangle -  \langle n T \nabla \cdot \mathbf{V} \rangle + \langle \eta \mathbf{j}, \mathbf{j} \rangle, \label{dT_dTdt}\\
&\left\langle \dd{H}{\mathbf{B}}, \pp{\mathbf{B}}{t} \right\rangle = - \left\langle \frac{\mathbf{B}}{\mu_0}, \nabla \times \mathbf{E} \right\rangle. \label{dB_dBdt}
\end{align}
\end{subequations}
\endgroup
The first term on the right-hand side of \eqref{dv_dvdt} vanishes due to the cross-product of $\mathbf{V}$ with itself. Further, upon applying integration by parts for the second and third term on the right-hand side of \eqref{dv_dvdt}, as well as the first terms on the right-hand sides of \eqref{dT_dTdt} and \eqref{dB_dBdt}, we arrive at
\begingroup
\addtolength{\jot}{2mm}
\begin{subequations}
\begin{align}
&\left\langle \dd{H}{n}, \pp{n}{t} \right\rangle = - \left\langle \frac{1}{2} m_i |\mathbf{V}|^2, \nabla \cdot (n \mathbf{V}) \right\rangle - \left\langle \frac{T}{\gamma - 1}, \nabla \cdot (n \mathbf{V}) \right\rangle,\\
&\left\langle \dd{H}{\mathbf{V}}, \pp{\mathbf{V}}{t} \right\rangle = \left\langle \frac{1}{2} m_i |\mathbf{V}|^2, \nabla \cdot (n \mathbf{V}) \right\rangle +  \langle n T \nabla \cdot \mathbf{V} \rangle + \left\langle \mathbf{V}, \mathbf{j} \times \mathbf{B} \right\rangle, \label{change_dHdV}\\
&\left\langle \dd{H}{T}, \pp{T}{t} \right\rangle = \left\langle \frac{T}{\gamma - 1}, \nabla \cdot (n \mathbf{V}) \right\rangle -  \langle n T \nabla \cdot \mathbf{V} \rangle + \langle \eta \mathbf{j}, \mathbf{j} \rangle, \label{change_dHdT} \\
&\left\langle \dd{H}{\mathbf{B}}, \pp{\mathbf{B}}{t} \right\rangle = - \left\langle \mathbf{j}, - \mathbf{V} \times \mathbf{B} + \eta \mathbf{j} \right\rangle, \label{ME_change}
\end{align}
\end{subequations}
\endgroup
where we used the boundary conditions, noting that we used the analogous electric field boundary condition $\mathbf{E} \times \mathbf{n} |_{\partial \Omega} = 0$ in place of the magnetic field one. Further, we used the definition of $\mathbf{j}$ in \eqref{dB_dBdt} after integrating by parts. Summing over the four terms as in \eqref{Hamiltonian_full}, we then find that all terms cancel, and the rate of change of $H$ in time equals zero. $\hfill \blacksquare$\\

The computations appearing in the above proof can also be used to uncover the system's transfers between the total kinetic, internal, and magnetic energies, which are given by
\begin{equation}
K\!E(n , \mathbf{V}) = \int_\Omega \frac{1}{2}m_i n |\mathbf{V}|^2dx, \hspace{1cm} I\!E(n , T) = \int_\Omega \frac{nT}{\gamma - 1}dx, \hspace{1cm} M\!E(\mathbf{B}) = \int_\Omega \frac{|\mathbf{B}|^2}{2\mu_0}dx.
\end{equation}
Using the variational derivative relations
\begin{equation}
\dd{H}{n} = \dd{K\!E}{n} + \dd{I\!E}{n}, \hspace{5mm} \dd{H}{T} = \dd{I\!E}{T}, \hspace{5mm} \dd{H}{\mathbf{V}} = \dd{K\!E}{\mathbf{V}}, \hspace{5mm} \dd{H}{\mathbf{B}} = \dd{M\!E}{\mathbf{B}},
\end{equation}
we then find
\begingroup
\addtolength{\jot}{2mm}
\begin{subequations} \label{H_transfers}
\begin{align}
&\frac{dK\!E}{dt} = - \langle \mathbf{V}, \nabla ( n T) \rangle  + \langle \mathbf{V},  \mathbf{j} \times \mathbf{B} \rangle, \\
&\frac{dI\!E}{dt} \;\;= + \langle \mathbf{V}, \nabla ( n T) \rangle \hspace{21mm}+ \langle \eta \mathbf{j}, \mathbf{j} \rangle \\
&\frac{dM\!E}{dt} = \hspace{23mm}- \langle \mathbf{V}, \mathbf{j} \times \mathbf{B} \rangle - \langle \eta \mathbf{j}, \mathbf{j} \rangle, \label{dt_ME}
\end{align}
\end{subequations}
\endgroup
and there is a kinetic-internal energy transfer based on the pressure gradient force, a kinetic-magnetic energy transfer based on the Lorentz force, and a magnetic-to-internal energy transfer based on resistivity.\\ \\
Next to energy conservation, we will consider two conserved quantities of the magnetic field equation. The first is given by the magnetic field's zero-divergence property.\\ \\
\textbf{Proposition 1b.} \textit{If a magnetic field governed by \eqref{B_eqn_cts} is initially divergence-free, it will remain so for all time.}\\ \vspace{-3mm} \\
\textbf{Proof.} Applying the divergence to the magnetic field evolution equation \eqref{B_eqn_cts}, we find
\begin{equation}
\pp{\nabla \cdot \mathbf{B}}{t} = - \nabla \cdot (\nabla \times \mathbf{E}) = 0,
\end{equation}
where we used the vector calculus identity $\nabla \cdot \big(\nabla \times (\cdot)\big) \equiv 0$. Since this shows that the rate of change of the magnetic field's divergence is zero, we find that the latter will always remain zero. $\hfill \blacksquare$\\

The second quantity is given by conservation of total magnetic helicity in the absence of resistivity and for appropriate boundary conditions. The total magnetic helicity is given by
\begin{equation}
H_M \coloneqq \langle \mathbf{A}, \mathbf{B} \rangle,
\end{equation}
where $\mathbf{A}$ denotes the magnetic vector potential, which is defined such that $\nabla \times \mathbf{A} = \mathbf{B}$. Note that $\mathbf{A}$ admits an evolution equation of the form
\begin{equation}
\pp{\mathbf{A}}{t} = - \mathbf{E} + \nabla \Phi, \label{A_eqn_cts}
\end{equation}
for any scalar potential field $\Phi$. Using the vector calculus identity $\nabla \times \big(\nabla (\cdot)\big) \equiv 0$, it is straightforward to verify that upon applying the curl operator to \eqref{A_eqn_cts}, the magnetic field equation \eqref{B_eqn_cts} is retrieved.\\ \\
\textbf{Proposition 1c.} \textit{Suppose \eqref{B_eqn_cts} is complemented by boundary conditions of the form $\mathbf{E} \times \mathbf{n}|_{\partial \Omega} = 0$, and further that $\eta = 0$. Then the total magnetic helicity $H_M$ is conserved.} \\ \vspace{-3mm} \\
\textbf{Proof.} Using the evolution equations \eqref{B_eqn_cts} and \eqref{A_eqn_cts}, we find
\begingroup
\addtolength{\jot}{2mm}
\begin{align}
\frac{d}{dt} H_M = \frac{d}{dt} \langle \mathbf{A}, \mathbf{B} \rangle = \left\langle \pp{\mathbf{A}}{t}, \mathbf{B} \right\rangle + \left\langle \mathbf{A},  \pp{\mathbf{B}}{t} \right\rangle = \langle -\mathbf{E} + \nabla \Phi, \mathbf{B} \rangle + \left\langle \mathbf{A},  - \nabla \times \mathbf{E} \right\rangle.
\end{align}
\endgroup
Applying integration by parts for the grad and curl terms, and using $\nabla \cdot \mathbf{B} = 0$, $\nabla \times \mathbf{A} = \mathbf{B}$, as well as $\mathbf{E} = - \mathbf{V} \times \mathbf{B} + \eta \mathbf{j}$, we then obtain
\begingroup
\addtolength{\jot}{2mm}
\begin{align}
\frac{d}{dt} \langle \mathbf{A}, \mathbf{B} \rangle = 2 \langle \mathbf{V} \times \mathbf{B}, \mathbf{B} \rangle - 2 \langle \eta \mathbf{j}, \mathbf{B} \rangle+ \int_{\partial \Omega} \Phi \mathbf{B} \cdot \mathbf{n} \; dS + \int_{\partial \Omega} \mathbf{A} \cdot (\mathbf{E} \times \mathbf{n}) dS. \label{Helicity_zero}
\end{align}
\endgroup
Finally, the first term vanishes due to the cross product of $\mathbf{B}$ with itself, and the remaining terms evaluate to zero given the proposition's assumptions (noting that the analogous boundary condition for $\mathbf{E} \times \mathbf{n} = 0$ is given by $\mathbf{B} \cdot \mathbf{n} = 0$).
$\hfill \blacksquare$ \\

In the above proof, we implicitly also showed that the rate of change of total helicity is invariant under the choice of scalar potential $\Phi$ only in the presence of appropriate boundary conditions. In cases where $\mathbf{B} \cdot \mathbf{n}$ is not equal to zero along some parts of the boundary $\partial \Omega$, a so-called relative magnetic helicity needs to be considered instead \cite{finn1985magnetic}. In the numerical results section below, we restrict ourselves to demonstrating helicity preservation in cases where considering $H_M$ -- rather than the relative helicity -- is sufficient.\\ \\
\textbf{Remark 2} In the presence of resistivity, the rate of change of total helicity reads
\begin{equation}
\frac{d}{dt} H_M = - 2 \langle \eta \mathbf{j}, \mathbf{B} \rangle. \label{dt_ME_k}
\end{equation}
In comparison, recalling \eqref{dt_ME}, we have a rate of change of total magnetic energy
\begin{equation}
\frac{dM\!E}{dt} = - \left\langle \mathbf{V}, \mathbf{j} \times \mathbf{B} \right\rangle - \left\langle \eta \mathbf{j}, \mathbf{j} \right\rangle.  \label{dt_HM_k}
\end{equation}
Since $\mathbf{j}$ is a first derivative of $\mathbf{B}$, we therefore find that current density features with a given wave number $k$ lead to a decay of total magnetic energy in the order of $O(k^2)$, and a decay of total magnetic helicity in the order of $O(k)$. In particular, for small scale or high wave number features, this behavior leads to a much faster decay in total magnetic energy than total magnetic helicity. In the limit of a strongly magnetized plasma with a tiny resistivity $\eta$ and where high $k$ modes dominate the plasma dynamics, the plasma tends to relax into a minimal magnetic energy state while holding the magnetic helicity conserved, a process known as magnetic relaxation or magnetic self-organization~\cite{taylor-prl-1974,yeates2020magnetohydrodynamic}. If a discretization of the resistive MHD equations does not satisfy an analogous discrete form of \eqref{dt_ME_k} and \eqref{dt_HM_k}, then the dynamics may relax to a different, unphysical minimum.

%
\subsection{Compatible finite element spaces}
In the following, we discretize equations \eqref{cty_eqn_full} - \eqref{B_eqn_cts} using the mixed compatible finite element method \cite{arnold2006finite, boffi2013mixed}, thereby ensuring that the vector calculus identities we used above hold discretely. For this purpose, we consider a discrete de Rham complex
\begin{displaymath}
    \xymatrix{
        H^1(\Omega) \ar[r]^{\nabla} \ar[d]^{\pi^0} & H(\text{{\normalfont curl}};\Omega) \ar[r]^{\nabla \times} \ar[d]^{\pi^1} & H(\text{{\normalfont div}};\Omega)\ar[d]^{\pi^2} \ar[r]^{\nabla \cdot} & L^2(\Omega) \ar[d]^{\pi^3} \\
      \mathbb{V}_0(\Omega) \ar[r]^{\nabla}  &\mathbb{V}_1(\Omega) \ar[r]^{\nabla \times} & \mathbb{V}_2(\Omega) \ar[r]^{\nabla \cdot}  & \mathbb{V}_3(\Omega)}
\end{displaymath}
for bounded projections $\pi_i$ and suitable finite element spaces $\mathbb{V}_i(\Omega)$, for $i \in (0, 1, 2, 3)$. The choice of finite element spaces depends on the underlying choice of mesh, and in the numerical results section below, we will consider a tetrahedral mesh with corresponding discrete de Rham complex
\begin{equation}
    \xymatrix{
      P_2(\Omega) \ar[r]^{\nabla}  &N1_2^e(\Omega) \ar[r]^{\nabla \times} & N1_2^f(\Omega) \ar[r]^{\nabla \cdot}  & dP_1(\Omega)}, \label{N1_complex}
\end{equation}
which is based on tetrahedral curl- and div-conforming Nedelec elements $N1$ of the first kind \cite{arnold2014periodic}\footnote{The superscripts $e$, $f$ denote elements based on moments integrated over edges and faces, respectively, which ensures that the elements are curl- and div-conforming, respectively.}, and where $P$ and $dP$ denote the standard continuous and discontinuous Galerkin spaces for tetrahedra, respectively. Subscripts correspond to the spaces' polynomial degree. Further, we will consider a hexahedral mesh with corresponding discrete de Rham complex
\begin{equation}
    \xymatrix{
      Q_2(\Omega) \ar[r]^{\nabla}  &Nc_2^e(\Omega) \ar[r]^{\nabla \times} & Nc_2^f(\Omega) \ar[r]^{\nabla \cdot}  & dQ_1(\Omega)}, \label{Q1_complex}
\end{equation}
which is based on the hexahedral curl- and div-conforming Nedelec elements $Nc$ of the first kind, and where $Q$ and $dQ$ denote the standard continuous and discontinuous Galerkin spaces for hexahedra. Given the complex, we set $n \in \mathbb{V}_3(\Omega)$, $\mathbf{V} \in \mathbb{V}_2(\Omega)$, and $T \in \mathbb{V}_0(\Omega)$. Further, we will derive discretizations for both a div-conforming magnetic field $\mathbf{B} \in \mathbb{V}_2(\Omega)$, and a curl-conforming magnetic field $\mathbf{B} \in \mathbb{V}_1(\Omega)$. Finally, in the presence of boundaries, we adjust the vector finite element spaces accordingly, that is
\begin{equation}
\mathbf{V} \in \mathring{\mathbb{V}}_2(\Omega) \subset \mathring{H}(\text{div}; \Omega) = \{\mathbf{w} \in H(\text{div}; \Omega) \; | \; \mathbf{w} \cdot \mathbf{n}_{|_{\partial \Omega}} = 0\}.
\end{equation}
For the magnetic and electric fields, we apply the normal and tangential boundary conditions \eqref{B_bc} either weakly or strongly according to the choice of function space. For div-conforming $\mathbf{B}$, we apply the normal boundary conditions on $\partial \Omega_1$ strongly, and therefore adjust the underlying div-conforming function space accordingly. Conversely, the tangential boundary conditions on $\partial \Omega_2$ are applied weakly in the computation of the curl-conforming current density $\mathbf{j} \in \mathbb{V}_1(\Omega)$. An analogous approach applies when $\mathbf{B}$ is set to be curl-conforming, with the tangential boundary conditions on $\partial \Omega_2$ applied strongly, and the ones on $\partial \Omega_1$ applied to $\mathbf{E}$ weakly. More details on this procedure are given in the next section.
%
\section{Space discretization} \label{sec_Space_discretization}
In order to discretize the ideal compressible MHD equations \eqref{MHD_cts} (with $\eta = 0$), we first consider the continuity, momentum, and temperature equations. Once we have discretized the latter equations and analyzed their energy conserving properties, we proceed to the magnetic field equation.
%
\subsection{Hydrodynamic system}
For the continuity and momentum equations, we follow \cite{natale2016compatible} up to the discretization of the pressure gradient and Lorentz force terms. For the pressure gradient, we find that $\nabla (n T)$ is ill-defined since $n \in \mathbb{V}_3$ is discontinuous across facets, and we therefore need to introduce an additional facet integral. Next, we consider the temperature equation. Recall that in the derivation of the conserved total energy, the pressure gradient term cancels with the divergence term in the temperature equation (see \eqref{change_dHdV} - \eqref{change_dHdT}). In order to achieve a discrete energy balance, we therefore need to discretize the latter term analogously to the pressure gradient term. Finally, we include an interior penalty term \cite{burman2006edge} in the temperature equation for transport stabilization purposes. Altogether, this leads to the following space discretization.\\ \\
\textbf{Definition 1.} The temperature based space discretization of the hydrodynamic equation set \eqref{cty_eqn_full} - \eqref{T_eqn_cts} solves for $(n, \mathbf{V}, T) \in \big(\mathbb{V}_3(\Omega), \mathring{\mathbb{V}}_2(\Omega), \mathbb{V}_0(\Omega)\big)$ such that
\begingroup
\allowdisplaybreaks
\addtolength{\jot}{2mm}
\begin{subequations} \label{nVT_discr}
\begin{align}
&\pp{n}{t} + \nabla \cdot \mathbf{F} = 0, \label{ideal_n_ec_T}\\
&\left\langle \mathbf{w}, \pp{\mathbf{V}}{t} \right\rangle + \left\langle \nabla \times \left( \mathbf{w} \times \mathbf{F}/n \right), \mathbf{V} \right\rangle - \int_\Gamma \{\!\{\mathbf{n} \times \left( \mathbf{w} \times \mathbf{F}/n \right)\}\!\} \cdot \tilde{\mathbf{V}} \; dS \nonumber \\
& \hspace{20mm} - \left\langle \frac{1}{2}|\mathbf{V}|^2, \nabla \cdot \mathbf{w} \right\rangle + \left\langle \mathbf{j}, \frac{\mathbf{w}}{m_i n} \times \mathbf{B}\right\rangle \nonumber\\ 
&\hspace{20mm} + \left\langle \frac{\mathbf{w}}{m_i n}, \nabla ( n T) \right\rangle - \int_\Gamma [\![ n T]\!] \left\{\frac{\mathbf{w}}{m_i n} \right\}dS = 0 &\forall \mathbf{w} \in \mathring{\mathbb{V}}_2(\Omega), \label{ideal_u_ec_upw}\\
& \left\langle \frac{ n \chi}{\gamma -1}, \pp{T}{t}  \right\rangle + \left\langle \frac{ \chi}{\gamma -1} \mathbf{F}, \nabla T \right\rangle + \int_\Gamma \frac{h_e^2 \kappa_T}{(\gamma - 1)\{n\}} [\![\mathbf{F} \cdot \nabla \chi]\!][\![\mathbf{F} \cdot \nabla T]\!]dS \nonumber \\
& \hspace{2cm} - \langle \mathbf{F}/n, \nabla (\chi  n T) \rangle + \int_\Gamma [\![ \chi  n T ]\!] \left\{ \mathbf{F}/n \right\}dS = 0 & \forall \chi \in \mathbb{V}_0(\Omega) \label{ideal_T_ec_ip},
\end{align}
\end{subequations}
\endgroup
for flux $\mathbf{F} \in \mathring{\mathbb{V}}_2(\Omega)$ defined according to
\begin{align}
\langle \mathbf{w}, \mathbf{F} - n\mathbf{V} \rangle = 0 && \forall \mathbf{w} \in \mathring{\mathbb{V}}_2(\Omega).
\end{align}
$\Gamma$ denotes the set of all interior facets of the mesh, and we applied jump, averaging, and upwind facet operations defined by
\begingroup
\addtolength{\jot}{2mm}
\begin{align}
&[\![\psi]\!] \coloneqq \psi^+ - \psi^-,  &\{\psi\} \coloneqq \frac{1}{2}\left(\psi^+ + \psi^- \right),\hspace{15mm}&\\
&[\![\mathbf{u}]\!] \coloneqq \mathbf{u}^+ \cdot \mathbf{n}^+ + \mathbf{u}^- \cdot \mathbf{u}^-, &\{\mathbf{u}\} \coloneqq \frac{1}{2} \left(\mathbf{u}^+ \cdot \mathbf{n}^+ - \mathbf{u}^- \cdot \mathbf{n}^-\right),& \\
&\{\!\{\mathbf{u}\}\!\} \coloneqq \mathbf{u}^+ \times \mathbf{n}^+ + \mathbf{u}^- \times \mathbf{n}^-
&\tilde{\mathbf{u}} \coloneqq 
\begin{cases}
\mathbf{u}^+ \;\;\; \text{if } \mathbf{u}^+ \cdot \mathbf{n}^+ < 0, \\
\mathbf{u}^- \;\;\; \text{otherwise},
\end{cases}\hspace{3mm} \label{def_upw}&
\end{align}
\endgroup
for any vector field $\mathbf{u}$ and scalar field $\psi$. $\mathbf{n}$ denotes the facet normal vector, and the two sides of each mesh facet are arbitrarily denoted by $+$ and $-$ (and hence $\mathbf{n}^+ = - \mathbf{n}^-$). Additionally, $h_e$ denotes the local facet area, and $\kappa_T$ is a free stabilization parameter.\begin{center}\vspace{-15mm}\end{center}\hrulefill\\

Note that no boundary terms occur in the discrete momentum equation, recalling that we either assume periodic domains, or free-slip boundary conditions for $\mathbf{V}$. Further, in view of energy conservation, we followed \cite{natale2016compatible} and introduced an auxiliary flux $\mathbf{F}$, and applied an upwind stabilized form of the velocity transport's curl part. Finally, note that we still need to specify $\mathbf{B}$ and $\mathbf{j}$ in the Lorentz force term. Before exploring the various possibilities to do so, we demonstrate the hydrodynamic system's discrete energy conservation property. In particular, we will find the interior penalty term included in the temperature equation \eqref{ideal_T_ec_ip} to be neutral with respect to the rate of change of total energy.\\ \\
\textbf{Proposition 2.} \textit{Ignoring the Lorentz force term, the hydrodynamic system of equations \eqref{nVT_discr} conserves the system's total hydrodynamic energy}
\begin{equation}
H(n, \mathbf{V}, T) =  \int_\Omega \left( \frac{1}{2} m_i n |\mathbf{V}|^2 + \frac{n T}{\gamma - 1} \right) dx.
\end{equation}
\textbf{Proof.} First, note that the above energy is discrete in the sense that it is a function of finite element fields. To prove that it is conserved, we proceed as in the continuous case, noting that after discretization, the Hamiltonian variations \eqref{H_variations} are given by their projections into the finite element spaces of the given fields with respect to which the Hamiltonian is varied. As before, we compute the terms arising in \eqref{Hamiltonian_full}, and find
\begingroup
\addtolength{\jot}{2mm}
\begin{align}
\frac{dH}{dt} &=\left\langle \dd{H}{n}, \pp{n}{t} \right\rangle+ \left\langle \dd{H}{\mathbf{V}}, \pp{\mathbf{V}}{t} \right\rangle + \left\langle \dd{H}{T}, \pp{T}{t} \right\rangle\\
&=\left\langle \frac{1}{2} m_i |\mathbf{V}|^2 +  \frac{T}{\gamma - 1}, \pp{n}{t} \right\rangle + \left\langle m_i \mathbf{F}, \pp{\mathbf{V}}{t} \right\rangle + \left\langle \frac{n}{\gamma-1}, \pp{T}{t} \right\rangle,
\end{align}
\endgroup
and we dropped the projections for the Hamiltonian variations in $n$ and $T$, noting that this is possible since they are paired in an inner product with fields of the same spaces. Using the discrete equations \eqref{ideal_u_ec_upw} and \eqref{ideal_T_ec_ip} with test functions
\begin{equation}
\mathbf{w} = m_i \mathbf{F} \in \mathring{\mathbb{V}}_2(\Omega), \hspace{1cm} \chi \equiv 1,
\end{equation}
and taking the inner product of the discrete continuity equation \eqref{ideal_n_ec_T} with
\begin{equation}
\frac{1}{2} m_i |\mathbf{V}|^2 +  \frac{T}{\gamma - 1},
\end{equation}
we find that similarly to the continuous case described in Proposiation 1a, all terms cancel except for the interior penalty term (and ignoring the Lorentz force term). In particular, all integrations by parts are already performed in the discrete equations, with the exception of the temperature equation's transport term. However, we find that for the latter term with $\chi \equiv 1$, integration by parts leads to
\begin{equation}
\left\langle \frac{1}{\gamma -1} \mathbf{F}, \nabla T \right\rangle = -\left\langle \frac{ T}{\gamma -1}, \nabla \cdot \mathbf{F} \right\rangle,
\end{equation}
and no additional facet or boundary integral terms occur due to the choice of flow field boundary conditions (which equally hold for $\mathbf{F} \in \mathring{\mathbb{V}}_1(\Omega)$), and since $T\!\cdot\!\mathbf{F}$ has continuous normal components across facets. Finally, the interior penalty term also evaluates to zero since for our choice of test function $\chi \equiv 1$, we find $\nabla 1 \equiv \mathbf{0}$. $\hfill\blacksquare$\\ \\
\textbf{Remark 3.} It is also possible to derive an energy conserving hydrodynamic system when considering $p$ instead of $T$ as the thermal field, using the pressure field evolution equation \eqref{p_eqn_cts} (with $\eta = 0$ for the ideal MHD case). In this case, due to the differential volume form-type transport of $p$, we consider the latter field in the discontinuous space $\mathbb{V}_3(\Omega)$. It is then convenient to consider a momentum equation weighted by the density $m_i n$
\begin{equation}
m_i n\left(\pp{\mathbf{V}}{t} + (\nabla \times \mathbf{V}) \times \mathbf{V} + \frac{1}{2} \nabla |\mathbf{V}|^2\right) + \nabla p - \mathbf{j} \times \mathbf{B} = 0,
\end{equation}
since this allows for a straightforward discretization of the pressure term, given by $-\langle \nabla \cdot \mathbf{w}, p \rangle$, where both $\nabla \cdot \mathbf{w}$ and $p$ are functions in $\mathbb{V}_3(\Omega)$. Using ideas described in \cite{wimmer2020density} for an energy conserving, upwind stabilized discretization of the density number and momentum equations, we then arrive at a hydrodynamic system of equations given by
\begingroup
\allowdisplaybreaks
\addtolength{\jot}{2mm}
\begin{subequations} \label{nVp_discr}
\begin{align}
&\left\langle \phi, \pp{n}{t} \right\rangle - \left\langle \nabla \phi, n \mathbf{V} \right\rangle + \int_\Gamma [\![\phi \mathbf{V}]\!] \tilde{n} \; dS = 0, \label{ideal_n_ec_p} \\
&\left\langle m_i n \; \mathbf{w}, \pp{\mathbf{V}}{t} \right\rangle + \left\langle \nabla \times \left( m_i n \; \mathbf{w} \times \mathbf{V} \right), \mathbf{V} \right\rangle - \int_\Gamma \{\!\{\mathbf{n} \times \left( m_i n\;\mathbf{w} \times \mathbf{V} \right)\}\!\} \cdot \tilde{\mathbf{V}} \; dS \nonumber \\
& \hspace{20mm} + \left\langle \nabla P, n\mathbf{w} \right\rangle - \int_\Gamma [\![P \mathbf{w}]\!] \tilde{n} \; dS + \left\langle \mathbf{j}, \mathbf{w} \times \mathbf{B}\right\rangle -\langle p, \nabla \cdot \mathbf{w} \rangle = 0 &\hspace{-1cm}\forall \mathbf{w} \in \mathring{\mathbb{V}}_2(\Omega), \\
&\left\langle \frac{\phi}{\gamma - 1}, \pp{p}{t}\right\rangle - \left\langle \frac{1}{\gamma - 1} \nabla \phi, p \mathbf{V} \right\rangle + \int_\Gamma [\![\phi\mathbf{V}]\!] \tilde{p} \; dS  +  \langle \phi p, \nabla \cdot \mathbf{V} \rangle = 0 & \forall \phi \in \mathbb{V}_3(\Omega) \label{ideal_p_ec_upw},
\end{align}
\end{subequations}
\endgroup
where the scalar upwind value $\tilde{n}$ is defined analogously to the vector version \eqref{def_upw}, and $P \in \mathbb{V}_3$ is defined by
\begin{align}
\left\langle \phi, P - \frac{m_i}{2}|\mathbf{V}|^2 \right\rangle && \forall \phi \in \mathbb{V}_3(\Omega).
\end{align}
Note that $P$ corresponds to the variational derivative of the total energy in $n$, where the latter is now given by
\begin{equation}
H(n, \mathbf{V}, p, \mathbf{B}) = \int_\Omega \left( \frac{1}{2} m_i n |\mathbf{V}|^2 + \frac{p}{\gamma-1} + \frac{|\mathbf{B}|^2}{2 \mu_0} \right) dx. \label{H_p}
\end{equation}
It can be shown that this system also conserves the total energy up to the Lorentz force term, where we proceed as for the temperature based discretization, and choose test functions given by
\begin{equation}
\phi = P \text{ (in \eqref{ideal_n_ec_p})}, \hspace{1cm} \mathbf{w} = \mathbf{V}, \hspace{1cm} \mathbf \phi \equiv 1 \text{ (in \eqref{ideal_p_ec_upw})}.
\end{equation}
Given the choice of test function for the pressure equation, we find that the upwind stabilized terms corresponding to $\nabla \cdot (\mathbf{pV})$ -- i.e. the second and third terms in \eqref{ideal_p_ec_upw} -- vanish since $\nabla 1 \equiv 0$ and $[\![\mathbf{V}]\!] = 0$. In particular, this implies that like the interior penalty term for the temperature based setup \eqref{ideal_T_ec_ip}, the upwind stabilized discretization for the pressure transport term is energetically neutral.\\ \\
Further, for the pressure based setup, we consider an upwind stabilized version of the continuity equation \eqref{ideal_n_ec_p}, while for the temperature based one, we consider a ``compatible'' one (with $\nabla \cdot \mathbf{F} \in \mathbb{V}_3(\Omega)$) given by \eqref{ideal_n_ec_T}. Unfortunately, it is difficult to implement an upwind stabilized continuity equation for the energy conserving temperature based setup. Unlike in the pressure based setup, when computing the rate of change of total energy according to \eqref{Hamiltonian_full}, we require terms arising from the continuity equation to cancel with terms of the temperature equation. Therefore, when including upwinding for the continuity equation, the temperature equation needs to be adjusted accordingly to achieve energy conservation, leading to an unfavorable transport discretization in the latter equation. When considering MHD scenarios with persistent sharp density gradients, the discretized continuity equation \eqref{ideal_n_ec_T} may lead to significant spurious oscillations, and in these cases it may therefore be advantageous to consider the pressure based setup \eqref{nVp_discr} instead. In practice, we have not found \eqref{ideal_n_ec_T} to lead to such issues, except in the presence of shock waves.\\ \\
Finally, the below discussion on the choice of magnetic field discretization works equally for either thermal field based setup, except for the choice of transport velocity in the discretization of the electric field \eqref{Ohms_law}. In order to ensure balanced energy transfers, when considering the temperature based setup \eqref{nVT_discr}, the magnetic field transport velocity is given by $\mathbf{F}/n$ -- which will be described further below -- while for the pressure based setup \eqref{nVp_discr}, it is given by $\mathbf{V}$.
%
\subsection{Magnetic field equation} \label{sec_B_eqn}
To complete the discretization, it remains to consider the magnetic field and its evolution equation. Next to energy conservation, we also need to ensure that the discrete magnetic field's zero-divergence is preserved. Further, we will discuss the resulting discrete magnetic helicity balance.\\ \\
For div-conforming $\mathbf{B} \in \mathbb{V}_2(\Omega)$ and curl-conforming $\mathbf{E}, \mathbf{j} \in \mathbb{V}_1(\Omega)$, the evolution equation and the magnetic field related quantities have been described e.g. in \cite{bochev2002matching}, and are reviewed in Definition 2 below.\\ \\
\textbf{Definition 2.} For the magnetic field's div-conforming space discretization, we seek $(\mathbf{B}, \mathbf{E}, \mathbf{j}) \in \big(\mathbb{V}_2(\Omega), \mathbb{V}_1(\Omega), \mathbb{V}_1(\Omega)\big)$ such that
\begingroup
\addtolength{\jot}{2mm}
\begin{subequations} \label{BHdiv_discr}
\begin{align}
&\pp{\mathbf{B}}{t} + \nabla \times \mathbf{E} = 0 \label{BHdiv_b_eqn}, \\
&\langle \mathbf{\Sigma}, \mathbf{E} + \mathbf{F}/n \times \mathbf{B} \rangle = 0 & \forall \mathbf{\Sigma} \in \mathbb{V}_1(\Omega), \label{BHdiv_E_eqn}\\
&\left\langle \mathbf{\Sigma}, \mathbf{j} \right\rangle = \left\langle \nabla \times \mathbf{\Sigma}, \frac{1}{\mu_0} \mathbf{B} \right\rangle + \int_{\partial \Omega_2} \mathbf{g}_2 \cdot \mathbf{\Sigma} \; dS = 0 & \forall \mathbf{\Sigma} \in \mathbb{V}_1(\Omega), \label{BHdiv_j_eqn}
\end{align}
\end{subequations}
\endgroup
where $\mathbf{g}_2$ is a function corresponding to the boundary condition for $\mathbf{H} \times \mathbf{n}|_{\partial \Omega_2}$. The boundary conditions on $\partial \Omega_1$ -- that is $\mathbf{B} \cdot \mathbf{n}|_{\partial \Omega_1} = g_1$ or an analogous version for $\mathbf{E} \times \mathbf{n}|_{\partial \Omega_1}$ -- are applied strongly to either $\mathbf{E} \in \mathbb{V}_1(\Omega)$ or $\mathbf{B} \in \mathbb{V}_2(\Omega)$.\begin{center}\vspace{-15mm}\end{center}\hrulefill\\

In the next series of propositions, we review the structure preserving properties of the div-conforming setup without transport stabilization. To demonstrate energy conservation, we recall that in the expansion \eqref{Hamiltonian_full} of the rate of change in time of total energy, the magnetic field transport term cancels with the Lorentz force term. The same cancellation still holds true after discretization.\\ \\
\textbf{Proposition 3a.} \textit{For  boundary conditions of the form $g_1 = 0$ on $\partial \Omega$, the space-discretized ideal compressible MHD equations \eqref{nVT_discr} + \eqref{BHdiv_discr} conserve the system's total energy \eqref{Hamiltonian}.}\\ \vspace{-3mm} \\
\textbf{Proof.} In Proposition 2, we have already shown that the system conserves energy up to terms related to the magnetic field. For these remaining terms, we find from the product rule \eqref{Hamiltonian_full},
\begin{align}
\frac{dH}{dt} =& \left\langle \dd{H}{\mathbf{V}}, \pp{\mathbf{V}}{t} \right\rangle + \left\langle \dd{H}{\mathbf{B}}, \pp{B}{t} \right\rangle + \dots = - \left\langle \mathbf{j}, \frac{m_i\mathbf{F}}{m_i n} \times \mathbf{B} \right\rangle + \left\langle \frac{\mathbf{B}}{\mu_0}, - \nabla \times \mathbf{E} \right\rangle + \dots, \label{Hamiltonian_full_BHdiv}
\end{align}
where we used the discrete momentum equation \eqref{ideal_u_ec_upw} with $\mathbf{w} = m_i \mathbf{F}$, and the discrete magnetic field equation \eqref{BHdiv_b_eqn} with $\mathbf{\Sigma} = \mathbf{B}/\mu_0$. Finally, we substitute $\mathbf{\Sigma} = \mathbf{E}$ in \eqref{BHdiv_j_eqn} (noting that $\partial \Omega_2 = \emptyset$), which leads to
\begin{equation}
\left\langle \frac{\mathbf{B}}{\mu_0}, - \nabla \times \mathbf{E} \right\rangle = -\langle \mathbf{E}, \mathbf{j} 
\rangle = -\langle - \mathbf{F}/n \times \mathbf{B}, \mathbf{j} \rangle =  \left\langle \mathbf{j}, \frac{m_i\mathbf{F}}{m_i n} \times \mathbf{B} \right\rangle.
\end{equation}
In particular, this term cancels with the term originating from the Lorentz force (first term on right-hand side of \eqref{Hamiltonian_full_BHdiv}), and we find that the total energy is preserved as required. $\hfill\blacksquare$\\

Since the proofs of Proposition 2 and 3a follow the steps of the non-discretized version described in Proposition 1a, we immediately find that the energy transfers \eqref{H_transfers} still hold after discretization. They are given by
\begingroup
\addtolength{\jot}{2mm}
\begin{subequations} \label{H_transfers_discr}
\begin{align}
&\frac{dK\!E}{dt} = + \left\langle \mathbf{F}/n, \nabla ( n T) \right\rangle - \int_\Gamma [\![ n T]\!] \left\{\mathbf{F}/n \right\}dS \;\;\;\;+ \langle \mathbf{F}/n,  \mathbf{j} \times \mathbf{B} \rangle, \\
&\frac{dI\!E}{dt} \;\;= - \left\langle \mathbf{F}/n, \nabla ( n T) \right\rangle + \int_\Gamma [\![ n T]\!] \left\{\mathbf{F}/n \right\}dS\\
&\frac{dM\!E}{dt} = \hspace{61mm} \;\;\;\;- \langle \mathbf{F}/n, \mathbf{j} \times \mathbf{B} \rangle,
\end{align}
\end{subequations}
\endgroup
and analogous consistent versions also hold true for the remaining magnetic field space discretizations described in this work.\\ \\
\textbf{Proposition 3b.} \textit{If the discrete, div-conforming magnetic field governed by \eqref{BHdiv_b_eqn} is initially divergence-free, it will remain so for all time.}\\ \vspace{-3mm} \\
\textbf{Proof.} The proof is identical to the non-discretized case given in Proposition 1b, noting that the vector calculus identity $\nabla \cdot \big(\nabla \times (\cdot)\big) \equiv 0$ still holds true after discretization, due to the choice of compatible finite element spaces. $\hfill \blacksquare$\\

Unlike total energy conservation and the zero-divergence condition, total magnetic helicity preservation does not follow naturally from the space discretization \eqref{BHdiv_discr}. As shown in \cite{hu2021helicity}, in order to additionally preserve total magnetic helicity, we need to replace the magnetic field occurring in the transport term by its projection into the curl-conforming space $\mathbb{V}_1(\Omega)$.\\ \\
\textbf{Proposition 3c.} \textit{Consider a space discretization of \eqref{B_eqn_cts} -\eqref{Ohms_law} given by
\begingroup
\addtolength{\jot}{2mm}
\begin{subequations} \label{BHdiv_discr_hel}
\begin{align}
&\pp{\mathbf{B}}{t} + \nabla \times \mathbf{E} = 0 \label{BHdiv_b_eqn_hel}, \\
&\langle \mathbf{\Sigma}, \mathbf{E} + \mathbf{F}/n \times P_{\mathbb{V}_1(\Omega)}(\mathbf{B}) \rangle = 0 & \forall \mathbf{\Sigma} \in \mathbb{V}_1(\Omega), \label{BHdiv_E_eqn_hel}\\
&\left\langle \mathbf{\Sigma}, \mathbf{j} \right\rangle = \left\langle \nabla \times \mathbf{\Sigma}, \frac{1}{\mu_0} \mathbf{B} \right\rangle = 0 & \forall \mathbf{\Sigma} \in \mathbb{V}_1(\Omega), \label{BHdiv_j_eqn_hel}
\end{align}
\end{subequations}
\endgroup
with boundary condition $g_1 = 0$ on $\partial \Omega$, and where $P_{\mathbb{V}_1(\Omega)}$ denotes projection into the curl-conforming space $\mathbb{V}_1(\Omega)$. Further, assume that there exists a magnetic vector potential $\mathbf{A} \in \mathbb{V}_1(\Omega)$ such that
\begin{equation}
\nabla \times \mathbf{A} = \mathbf{B}, \hspace{1cm} \pp{\mathbf{A}}{t} = -\mathbf{E} + \nabla \Phi,
\end{equation}
for some scalar potential $\Phi \in \mathbb{V}_0(\Omega)$. Then the system's total magnetic helicity $H_M = \langle \mathbf{A}, \mathbf{B} \rangle$ is conserved.}\\ \vspace{-3mm} \\
\textbf{Proof.} First, we note that the above definition of $\mathbf{A}$ is well-defined. Since the latter field is curl-conforming, and $\mathbf{B}$ is div-conforming, the equation $\nabla \times \mathbf{A} = \mathbf{B}$ holds strongly. Further, the time evolution equation of $\mathbf{A}$ is consistent with \eqref{BHdiv_b_eqn_hel}, which follows from applying the curl and using the vector calculus identity $\nabla \times \big(\nabla (\cdot)\big) \equiv 0$. For the remaining proof, we follow the one for Proposition 1c and obtain
\begingroup
\addtolength{\jot}{2mm}
\begin{align}
\frac{d}{dt} \langle \mathbf{A}, \mathbf{B} \rangle = \left\langle \pp{\mathbf{A}}{t}, \mathbf{B} \right\rangle + \left\langle \mathbf{A},  \pp{\mathbf{B}}{t} \right\rangle  = \langle -\mathbf{E} + \nabla \Phi, \mathbf{B} \rangle + \left\langle \mathbf{A},  - \nabla \times \mathbf{E} \right\rangle,
\end{align}
\endgroup
noting that both the $\mathbf{A}$ and $\mathbf{B}$ evolution equations hold strongly. Proceeding as before, we next apply integration by parts for the grad and curl terms. The expressions $\Phi \mathbf{B} \cdot \mathbf{n}$ and $(\mathbf{A} \times \mathbf{E}) \cdot \mathbf{n}$ are continuous across facets given the choice of function spaces, and therefore no additional facet integrals appear after integrating by parts. Using $\nabla \cdot \mathbf{B} = 0$ and $\nabla \times \mathbf{A} = \mathbf{B}$ as well as the boundary conditions, we then obtain
\begingroup
\addtolength{\jot}{2mm}
\begin{align}
\frac{d}{dt} \langle \mathbf{A}, \mathbf{B} \rangle = -2 \langle \mathbf{B}, \mathbf{E} \rangle = -2 \left\langle P_{\mathbb{V}_1(\Omega)}(\mathbf{B}), \mathbf{E} \right\rangle = -2 \left\langle P_{\mathbb{V}_1(\Omega)}(\mathbf{B}), -\mathbf{F}/n \times P_{\mathbb{V}_1(\Omega)}(\mathbf{B}) \right\rangle = 0,
\end{align}
\endgroup
where we made use of projection properties (with $\mathbf{E} \in \mathbb{V}_1(\Omega)$) as well as the vanishing cross-product $P_{\mathbb{V}_1(\Omega)}(\mathbf{B}) \times P_{\mathbb{V}_1(\Omega)}(\mathbf{B})$. $\hfill \blacksquare$\\

Next to the above div-conforming setup for $\mathbf{B}$, it is also possible to consider the latter field in the curl-conforming space, following e.g. \cite{bossavit1982mixed, schotzau2004mixed}. The below definition reviews the corresponding setup.\\ \\
\textbf{Definition 3.} For the magnetic field's curl-conforming space discretization, we seek $\mathbf{B} \in \mathbb{V}_1(\Omega)$ such that
\begin{align}
\left\langle \mathbf{\Sigma}, \pp{\mathbf{B}}{t} \right\rangle + \left\langle \nabla \times \mathbf{\Sigma}, - \mathbf{F}/n \times \mathbf{B} \right\rangle + \int_{\partial \Omega_1} \mathbf{\Sigma} \cdot \mathbf{h}_1 \; dS  = 0 && \forall \mathbf{\Sigma} \in \mathbb{V}_1(\Omega), \label{BHcurl_discr}
\end{align}
where $\mathbf{h}_1$ corresponds to the boundary condition for $\mathbf{E} \times \mathbf{n}|_{\partial \Omega_1}$. Further, the boundary condition on $\partial \Omega_2$ -- that is a function $\mathbf{g}_2$ corresponding to $\mathbf{H} \times \mathbf{n}|_{\partial \Omega_2}$ -- is applied strongly to $\mathbf{B} \in \mathbb{V}_1(\Omega)$.\begin{center}\vspace{-12mm}\end{center}\hrulefill\\

The -- in this case div-conforming -- electric field $\mathbf{E} \in \mathbb{V}_2(\Omega)$ is defined by
\begin{align}
&\langle \mathbf{w}, \mathbf{E} + \mathbf{F}/n \times \mathbf{B} \rangle = 0 & \forall \mathbf{\Sigma} \in \mathbb{V}_2(\Omega). \label{BHcurl_E_eqn}
\end{align}
However, it does not need to be computed explicitly as a projection for the magnetic field equation \eqref{BHcurl_discr}, since it is paired with the div-conforming field $\nabla \times \mathbf{\Sigma} \in \mathbb{V}_2(\Omega)$. Equally, the density current $\mathbf{j} \equiv 1/\mu_0 \nabla \times \mathbf{B} \in \mathbb{V}_2(\Omega)$ does not need to be computed as a separate projection, since the curl-operator maps the magnetic field strongly into the div-conforming space.\\ \\
As for the div-conforming magnetic field setup, we next move on to reviewing the structure preserving properties of the curl-conforming setup without transport stabilization.\\ \\
\textbf{Proposition 4a.} \textit{For  boundary conditions of the form $\mathbf{h}_1 = 0$ on $\partial \Omega$, the space-discretized ideal compressible MHD equations \eqref{nVT_discr} + \eqref{BHcurl_discr} conserve the system's total energy \eqref{Hamiltonian}.}\\ \vspace{-3mm} \\
\textbf{Proof.} As before, from Proposition 2, we have already shown that the system conserves energy up to terms related to the magnetic field. For these remaining terms, we find from the product rule \eqref{Hamiltonian_full},
\begin{align}
\frac{dH}{dt} =& \left\langle \dd{H}{\mathbf{V}}, \pp{\mathbf{V}}{t} \right\rangle + \left\langle \dd{H}{\mathbf{B}}, \pp{\mathbf{B}}{t} \right\rangle + \dots = - \left\langle \mathbf{j}, \frac{m_i\mathbf{F}}{m_i n} \times \mathbf{B} \right\rangle + \left\langle \nabla \times \frac{\mathbf{B}}{\mu_0}, \mathbf{F}/n \times \mathbf{B} \right\rangle + \dots, \label{Hamiltonian_full_BHcurl}
\end{align}
where we used the discrete momentum equation \eqref{ideal_u_ec_upw} with $\mathbf{w} = m_i \mathbf{F}$, and the discrete magnetic field equation \eqref{BHcurl_discr} with $\mathbf{\Sigma} = \mathbf{B}/\mu_0$. Note that in the momentum equation, we have $\mathbf{j} = 1/\mu_0 \nabla \times \mathbf{B} \in \mathbb{V}_2(\Omega)$. Making the same substitution for $\mathbf{j}$ in the last term of \eqref{Hamiltonian_full_BHcurl}, we obtain the required cancellation. $\hfill\blacksquare$\\ \\
\textbf{Proposition 4b.} \textit{Consider the weak magnetic divergence $\delta_B \in \mathbb{V}_0(\Omega)$, defined according to
\begin{align}
&\left\langle \chi, \delta_B \right\rangle = - \left\langle \nabla \chi, \mathbf{B} \right\rangle + \int_{\partial\Omega} \chi g_1 \; dS & \forall \chi \in \mathbb{V}_0(\Omega), \label{def_weak_B_div}
\end{align}
where $\mathbf{B}$ is governed by \eqref{BHcurl_discr}, and $g_1$ is defined on $\partial \Omega$ such that
\begin{equation}
\pp{g_1}{t} - \nabla \cdot \mathbf{h}_1 = 0. \label{Bn_Exn_cc}
\end{equation}
Then if the initial weak magnetic divergence is zero, it will remain so for all time.}\\ \vspace{-3mm} \\
\textbf{Proof.} First, we note that the above definition of $\delta_B$ is consistent with the non-discretized magnetic field and its divergence. To see this, we take the inner product of the non-discretized initial zero-divergence relation $\nabla \cdot \mathbf{B} = 0$ with a test function $\chi \in V_0(\Omega)$, and apply integration by parts to arrive at
\begin{align}
\left\langle \chi, \nabla \cdot \mathbf{B} \right\rangle = - \left\langle \nabla \chi, \mathbf{B} \right\rangle + \int_{\partial\Omega} \chi (\mathbf{B} \cdot \mathbf{n}) \; dS && \forall \chi \in \mathbb{V}_0(\Omega). \label{def_weak_non_discr_B_div}
\end{align}
Finally, we note that \eqref{Bn_Exn_cc} corresponds to the strong boundary compatibility equation of the form 
\begin{equation}
\pp{\mathbf{B} \cdot \mathbf{n}}{t} - \nabla \cdot (\mathbf{E} \times \mathbf{n}) = 0,
\end{equation}
which can be derived from the non-discretized magnetic field equation (and noting that $\nabla \times \mathbf{n}$ = 0). Therefore, since $\mathbf{h}_1$ corresponds to $\mathbf{E} \times \mathbf{n}$ by definition, we find that $g_1$ as defined in \eqref{Bn_Exn_cc} does to $\mathbf{B} \cdot \mathbf{n}$, implying that \eqref{def_weak_B_div} is consistent with the non-discretized, weak divergence-free condition \eqref{def_weak_non_discr_B_div}.\\ \\
To prove the proposition's claim, we take the time derivative of \eqref{def_weak_B_div}, and use \eqref{BHcurl_discr} to arrive at
\begin{align}
\left\langle \chi, \pp{\delta_B}{t} \right\rangle = \left\langle \nabla \times\nabla \chi, - \mathbf{F}/n \times \mathbf{B} \right\rangle + \int_{\partial \Omega_1} \nabla \chi \cdot \mathbf{h}_1 \; dS + \int_{\partial\Omega} \chi \pp{g_1}{t} \; dS && \forall \chi \in \mathbb{V}_0(\Omega), \label{def_weak_B_div_dt_pre}
\end{align}
where we could substitute for \eqref{BHcurl_discr} since $\nabla \chi \in \mathbb{V}_1(\Omega)$. Using the vector calculus identity $\nabla \times \nabla (\cdot) \equiv 0$, we find that the first term on the right-hand side vanishes. Finally, using integration by parts for the boundary integral involving $\mathbf{h}_1$ and the defining relation \eqref{Bn_Exn_cc}, we find that 
\begin{align}
\left\langle \chi, \pp{\delta_B}{t} \right\rangle = 0 && \forall \chi \in \mathbb{V}_0(\Omega), \label{def_weak_B_div_dt}
\end{align}
and the result follows upon setting $\chi \equiv \delta_B$. $\hfill \blacksquare$\\

Note that when integrating by parts for $\mathbf{h}_1$ in \eqref{def_weak_B_div_dt_pre}, a sufficient amount of regularity for $\mathbf{h}_1$ has implicitly been assumed. If $\mathbf{h}_1$ is not sufficiently regular, a weak form of \eqref{Bn_Exn_cc} has to be used instead in order to compute $g_1$ used in \eqref{def_weak_B_div}.\\ \\
Further, note that with \eqref{def_weak_B_div_dt} we strictly speaking showed something stronger than total magnetic field divergence preservation. In particular, the latter equation implies that the divergence is also weakly preserved in a local sense due to the inner product of $\pp{\delta_B}{t}$ with any test function $\chi$. We would therefore expect the curl-conforming setup to work as well as the div-conforming one in terms of controlling the magnetic field's divergence, up to jumps in the normal component of $\mathbf{B}$ across facets. The latter need not be zero for curl-conforming $\mathbf{B}$, and may grow over time. In our numerical results, we found that the magnetic field transport stabilization methods to be introduced below also effectively controlled and mitigated the magnitude of such jumps.\\ \\
It remains to discuss total helicity preservation for the curl-conforming setup \eqref{BHcurl_discr}. For the div-conforming case, we showed that a total helicity based on a curl-conforming magnetic vector potential $\mathbf{A}$ can be preserved if the underlying magnetic field evolution equation is adjusted appropriately. Analogously, one would expect that for the curl-conforming case, a similar approach might be possible based on a div-conforming magnetic vector potential, from which the magnetic field can be derived using a weak curl operation. However, we could not derive an adjusted, total helicity preserving setup akin to the one in Proposition 3c for this supposed ``mirror setup''. Instead, we had to revert to a setup in which $\mathbf{E}$, $\mathbf{j}$, and $\mathbf{A}$ are also curl-conforming and are defined as in Proposiction 3c. While this approach does conserve magnetic helicity, we found it to perform badly when compared to the original curl-conforming setup \eqref{BHcurl_discr}, and omit it here. In future work, it may be of interest to study further why the latter curl-conforming setup with a div-conforming magnetic vector potential fails to preserve magnetic helicity -- i.e. where exactly the ``mirror setup'' argument fails -- for instance by considering a differential form framework as described e.g. in \cite{pagliantini2016computational}.
%
\subsubsection{Magnetic field transport stabilization} \label{sec_B_eqn_stab}
Having reviewed the structural properties of the div- and curl-conforming space discretizations for the magnetic field evolution equation, we next consider possible ways of including transport stabilization, which forms the core part of this work. While energy conservation and helicity preservation are only infrequently considered in $\mathbf{B}$ field discretizations, the zero-divergence property is always considered in an exact or approximate way in order to ensure numerical stability. When devising magnetic field transport stabilization methods, our primary concern is therefore to ensure that the latter property still holds exactly in our compatible finite element discretization. In the above, we reviewed that the key property was the discrete vector calculus identities $\nabla \cdot \big(\nabla \times (\cdot) \big) \equiv 0$ and $\nabla \times \nabla (\cdot) \equiv 0$. In particular, these hold independently of the exact form of the electric field $\mathbf{E}$, which is defined by \eqref{BHdiv_E_eqn} and \eqref{BHdiv_E_eqn_hel} for the div-conforming and helicity preserving, div-conforming setups, respectively, and by \eqref{BHcurl_E_eqn} for the curl-conforming setup. In particular, we will therefore aim to adjust the definition of $\mathbf{E}$, thereby not compromising the aforementioned discrete vector calculus identities that ensure exact zero-divergence after discretization. The resulting stabilized discretizations will then be equal to the ones presented in Definitions 2 and 3 (as well as the helicity preserving one defined in Proposition 3c if applicable), up to the definition of $\mathbf{E}$.\\

\textbf{SUPG based stabilization}\\ \\
For the first type of stabilization, we note that the magnetic field transport equation
\begin{equation}
\pp{\mathbf{B}}{t} + \nabla \times (\mathbf{B} \times \mathbf{F}/n) = 0, \label{strong_B_res}
\end{equation}
is of the same transport form as the vorticity equation (up to pressure and Lorentz force related terms), which can be obtained by taking the curl of the momentum equation \eqref{momentum_eqn}. In particular, vorticity stabilization methods derived in the literature in the context of compatible finite elements are immediately applicable here. This includes the Anticipated Potential Vorticity Method (APVM; \cite{sadourny1985parameterization}), and its extensions to SUPG; see \cite{bauer2018energy} for the 2D case, and an extension to 3D in \cite{wimmer2020energy}. Here, we follow the latter work and define the magnetic field equation's strong residual according to \eqref{strong_B_res}, that is
\begin{align}
\mathbf{B}_{res} = \pp{\mathbf{B}}{t} + \nabla \times (\mathbf{B} \times \mathbf{F}/n).
\end{align}
While the above equals zero before discretization, this need not be the case after discretization (with the curl considered cell-wise).\\ \\
\textbf{Definition 4.} The SUPG-adjusted electric field is given by
\begin{equation}
\mathbf{E}_{SU\!PG} \coloneqq P_{\mathbb{V}_i(\Omega)} \Big(- \mathbf{F}/n \times \big(\mathbf{B} - \tau \mathbf{B}_{res} \big)\Big), \label{E_SUPG}
\end{equation}
where $P_{\mathbb{V}_i(\Omega)}$ denotes the $L^2$-projection into $\mathbb{V}_i(\Omega)$, and $i=1$ for a curl-conforming electric field, and $i=2$ for a div-conforming one. Further, the stabilization parameter $\tau$ is given by
\begin{equation}
\tau = \left(\left(\frac{2 \lambda}{\Delta t}\right)^2 + \left(\frac{2|\mathbf{V}|}{h_c}\right)^2 \right)^{\frac{1}{2}},
\end{equation}
for time step $\Delta t$, local cell size $h_c$, and a tuning parameter $\lambda$ to be specified.
\begin{center}\vspace{-12mm}\end{center}\hrulefill\\

For the div-conforming discretization \eqref{BHdiv_discr}, we then have $i=1$ and \eqref{E_SUPG} replaces the original definition of $\mathbf{E}$, given by \eqref{BHdiv_E_eqn}. For the curl-conforming discretization \eqref{BHcurl_discr}, we instead have $i=2$ and \eqref{E_SUPG} replaces \eqref{BHcurl_E_eqn}. As before, for the curl-conforming case, the projection to compute $\mathbf{E}$ is not required due to the pairing of the latter with $\nabla \times \mathbf{\Sigma} \in \mathbb{V}_2(\Omega)$ in the magnetic field evolution equation.\\ \\
To illustrate the above formulation's potential stabilizing effect, we make the following observation using the curl-conforming discretization \eqref{BHcurl_discr} with \eqref{E_SUPG}. Restricting the equation to a 2D manifold $\tilde{\Omega}$ with $\mathbf{B}$ aligned perpendicular to the latter, we find that $\mathbf{B}$ restricted to $\tilde{\Omega}$ can be identified as a scalar field $B$. On $\tilde{\Omega}$, the underlying curl-conforming space $\mathbb{V}_1(\Omega)$ degenerates to the continuous Galerkin space $\mathbb{V}_0(\tilde{\Omega})$\footnote{Strictly speaking for this to be the case, given a finite element mesh $\mathcal{T}$ used for $\Omega$, we need to assume $\tilde{\Omega}$ to cut orthogonally and half-way through cell edges in $\mathcal{T}$. Further, we assume a lowest order curl-conforming space $\mathbb{V}_1(\Omega)$ with moments identified by point evaluations. Altogether, this ensures that the degrees of freedom of $\mathbb{V}_1(\Omega)$ are evaluated on, and point orthogonally to $\tilde{\Omega}$.}. The resulting restricted magnetic field equation is given by
\begin{align}
\left\langle \Sigma, \pp{B}{t} \right\rangle + \left\langle -\nabla^\perp \Sigma, \tilde{\mathbf{V}}^\perp (B - \tau B_{res}) \right\rangle = 0 && \forall \Sigma \in \mathbb{V}_0(\tilde{\Omega}),
\end{align}
where $(a, b)^\perp = (-b, a)$, and $\tilde{\mathbf{V}} \in \mathbb{V}_1(\tilde{\Omega})$ denotes the (projected) restriction of $\mathbf{F}/n$ onto $\tilde{\Omega}$ (and we also apply this restriction to $\mathbf{F}/n$ in $B_{res}$). Finally, following the vorticity-based argument in \cite{bauer2018energy}, the above can be reformulated to
\begin{align}
\left\langle \Sigma + \tau \tilde{\mathbf{V}} \cdot \nabla \Sigma, \pp{B}{t} + \nabla \cdot (\tilde{\mathbf{V}} B) \right\rangle = 0 && \forall \Sigma \in \mathbb{V}_0(\tilde{\Omega}). \label{B_curl_2D_SUPG}
\end{align}
The above equation corresponds to the 2D scalar magnetic field equation with a standard SUPG modification for the test functions $\Sigma$, thereby motivating our naming convention for the modified electric field $\mathbf{E}_{SU\!PG}$. In particular, this implies that at least along planes perpendicular to the magnetic field, the residual based SUPG modification \eqref{E_SUPG} will have a stabilizing effect given by a non-positive definite term of the form
\begin{align}
\frac{1}{2}\frac{d}{dt}\|B\|^2_2 = \dots - \|\sqrt{\tau} \tilde{\mathbf{V}} \cdot \nabla B\|^2_2. \label{SUPG_neg_def}
\end{align}
For the div-conforming discretization \eqref{BHdiv_discr}, the derivation follows less naturally, since in this case, the div-conforming space $\mathbb{V}_2(\Omega)$ degenerates to the discontinuous Galerkin space $\mathbb{V}_3(\tilde{\Omega})$. The projection for $\mathbf{E}_{SU\!PG}$ can then not be dropped, thereby prohibiting a reformulation which leads to a standard SUPG form akin to \eqref{B_curl_2D_SUPG}. This in turn is required to demonstrate the SUPG method's underlying stabilizing effect.\\ \\
From a structure preserving point of view, given the SUPG-stabilized electric field $\mathbf{E}_{SU\!PG}$, we next aim to restore energy conservation for both the div- and curl-conforming discretizations of $\mathbf{B}$ by modifying the overall discretization for the MHD equations appropriately. For this purpose, we reconsider the proofs of Propositions 3a and 4a. When using $\mathbf{E}_{SU\!PG}$, we obtain
\begin{equation}
\frac{dH}{dt} = \dots + \left\langle \dd{H}{\mathbf{B}}, \pp{\mathbf{B}}{t} \right\rangle = \dots + \left\langle \mathbf{j}, \mathbf{F}/n \times \mathbf{B} - \tau \mathbf{B}_{res} \right\rangle,
\end{equation}
for both the div- and curl-conforming discretizations of $\mathbf{B}$. In the aforementioned proofs, we recall that this term cancels with the one arising from the Lorentz force. To restore energy conservation, we therefore simply replace $\mathbf{B}$ in the latter force term by $\mathbf{B} - \tau \mathbf{B}_{res}$. The corresponding term in the momentum equation is then given by
\begin{align}
\left\langle \mathbf{w}, \pp{\mathbf{V}}{t} \right\rangle + \dots + \left\langle \mathbf{j}, \frac{\mathbf{w}}{m_i n} \times \big(\mathbf{B} - \tau \mathbf{B}_{res}\big)\right\rangle = 0 && \forall \mathbf{w} \in \mathring{\mathbb{V}}_2(\Omega), \label{Lorentz_modified}
\end{align}
and the proof for energy conservation is identical to Propositions 3a and 4a. Altogether, the modifications lead to a revised discrete energy budged, which instead of \eqref{H_transfers_discr} now reads
\begingroup
\addtolength{\jot}{2mm}
\begin{subequations}
\begin{align}
&\frac{dK\!E}{dt} = \dots + \langle \mathbf{F}/n,  \mathbf{j} \times \big(\mathbf{B} - \tau \mathbf{B}_{res}\big) \rangle, \\
&\frac{dM\!E}{dt} = \hspace{5.3mm}- \langle \mathbf{F}/n, \mathbf{j} \times \big(\mathbf{B} - \tau \mathbf{B}_{res}\big) \rangle.
\end{align}
\end{subequations}
\endgroup
Note that \eqref{Lorentz_modified} is only one possible choice to restore energy conservation given the SUPG-type stabilization in $\mathbf{E}_{SU\!PG}$. Given the stabilization's dissipative nature, one could instead add a corresponding term in the temperature equation to transfer the loss of magnetic energy due to the SUPG stabilization into internal energy, rather than kinetic energy. However, here we chose the modification \eqref{Lorentz_modified}, since the stabilization term directly appears in the magnetic field transport term, whose energetic counterpart is given by the Lorentz force term, rather than a term affecting the internal energy.\\ \\
Finally, we note that this type of stabilization is not compatible with helicity conservation, whose proof \eqref{Helicity_zero} relies on the cancellation of
\begin{equation}
\langle \mathbf{B} \times \mathbf{F}/n, \mathbf{B} \rangle,
\end{equation}
up to additional projections introduced after discretization. However, with the above SUPG modification of the transported magnetic field, the term reads
\begin{equation}
\langle \big(\mathbf{B} - \tau \mathbf{B}_{res}\big) \times \mathbf{F}/n, \mathbf{B} \rangle,
\end{equation}
(again up to projections introduced after discretization), and $\tau \mathbf{B}_{res} \times \mathbf{B}$ not be equal to zero.\\

\textbf{Sub-grid resistivity based stabilization}\\ \\
For the second type of stabilization, we recall the MHD equation's resistive version, which includes Ohmic heating in the thermal field equation \eqref{T_eqn_cts}, and a resistive contribution in the electric field \eqref{Ohms_law}. While the ideal MHD equations do not contain such resistive effects, their dynamics as seen through a finite resolution filter may. In particular, magnetic field lines exhibiting features smaller than the filter-width will appear to be reconnecting from the filtered point of view, which is a process associated with resistivity. For a discretized system in a fixed Eulerian grid, one may therefore aim to emulate this effect using some form of sub-grid scale resistivity. To achieve this, we replace the resistive contribution $\eta \mathbf{j}$ by an interior penalty \cite{burman2006edge} version thereof. Since this occurs within the electric field $\mathbf{E}$, as before this will not affect the discrete magnetic field's zero-divergence property.\\ \\
\textbf{Definition 5.} The electric field $\mathbf{E}_\eta$ adjusted to contain a sub-grid scale resistive effect is given by
\begin{align}
\left\langle \bs{\omega}, \mathbf{E}_{\eta} \right\rangle = \left\langle \bs{\omega},  - \mathbf{F}/n \times \mathbf{B} \right\rangle + \int_\Gamma h_e^2 \kappa_B [\![\bs{\omega}]\!]' \cdot [\![\mathbf{j}]\!]' dS && \forall \bs{\omega} \in \mathbb{V}_i(\Omega), \label{E_eta}
\end{align}
where $k_B$ is a suitable stabilization parameter, and $i=1$ for a curl-conforming electric field, and $i=2$ for a div-conforming one. Further, the vector jump $[\![ \cdot ]\!]'$ is defined by
\begin{equation}
[\![ \mathbf{u} ]\!]' = \mathbf{u}^+ - \mathbf{u}^-
\end{equation}
for any vector field $\mathbf{u}$. \begin{center}\vspace{-12mm}\end{center}\hrulefill\\
When considering the curl-conforming discretization \eqref{BHcurl_discr} for $\mathbf{B}$ together with $\mathbf{E}_\eta$, we find upon setting the test function to $\mathbf{\Sigma} = \mathbf{B}$,
\begin{align}
\frac{d}{dt} \int_\Omega \frac{|\mathbf{B}|^2}{2\mu_0} \; dx = \left\langle \frac{\mathbf{B}}{\mu_0}, \pp{\mathbf{B}}{t} \right\rangle = - \left\langle \frac{1}{\mu_0} \nabla \times \mathbf{B}, \mathbf{E}_\eta \right\rangle = -\langle \mathbf{j}, \mathbf{E}_\eta \rangle = \left\langle \mathbf{j},  \mathbf{F}/n \times \mathbf{B} \right\rangle - \int_\Gamma h_e^2 \kappa_B |[\![\mathbf{j}]\!]'|^2 dS.
\end{align}
In other words, the newly introduced non-positive definite facet integral term has a stabilizing effect by dissipating magnetic energy at the grid-scale, which is achieved through penalizing the cross-facet jumps of $\mathbf{j}$. Note that unlike a resistive contribution $\eta \mathbf{j}$ to the electric field, the interior penalty term is consistent in the sense that the jump term $[\![\mathbf{j}]\!]'$ vanishes for a strong solution with a continuous density current $\mathbf{j}$. Similarly, for the div-conforming discretization \eqref{BHdiv_discr} for $\mathbf{B}$ together with $\mathbf{E}_\eta$, we find
\begingroup
\allowdisplaybreaks
\begin{align}
\frac{d}{dt} \int_\Omega \frac{|\mathbf{B}|^2}{2\mu_0} \; dx &= \left\langle \frac{\mathbf{B}}{\mu_0}, \pp{\mathbf{B}}{t} \right\rangle = - \left\langle \frac{\mathbf{B}}{\mu_0}, \nabla \times \mathbf{E}_\eta \right\rangle = - \langle \mathbf{j}, \mathbf{E}_\eta \rangle + \int_{\partial \Omega_2} \mathbf{g}_2 \cdot \mathbf{E}_\eta \; dS \label{Bdiv_E_eta1}\\
&= \left\langle \mathbf{j},  \mathbf{F}/n \times \mathbf{B} \right\rangle - \int_\Gamma h_e^2 \kappa_B |[\![\mathbf{j}]\!]'|^2 dS + \int_{\partial \Omega_2} \mathbf{g}_2 \cdot \mathbf{E}_\eta \; dS, \label{Bdiv_E_eta2}
\end{align}
\endgroup
where we used the weak definition \eqref{BHdiv_j_eqn} of $\mathbf{j}$. As for the curl-conforming discretization for $\mathbf{B}$, we find that the facet integral term has a stabilizing effect.\\ \\
From a structure preserving point of view, we need to account for the amount of magnetic energy dissipated by the non-positive definite interior penalty term at the grid-scale. In line with the above filter-based argument, this has been done e.g. for kinetic energy dissipation in the ocean modeling community by including a backscatter term to re-inject energy at resolved scales \cite{jansen2014parameterizing, jansen2015energy}. However, for simplicity, here we follow the physical energy balance occurring in the resistive MHD equations, and transfer magnetic energy lost by resistivity into internal energy via Ohmic heating (last term in \eqref{T_eqn_cts}). The corresponding term in the discrete temperature equation\footnote{For the discrete pressure equation \eqref{ideal_p_ec_upw}, the corresponding Ohmic heating term is identical up to the test function $\phi$.} \eqref{ideal_T_ec_ip} is then given on the right-hand side by
\begin{align}
\left\langle \frac{ n \chi}{\gamma -1}, \pp{T}{t}  \right\rangle + \dots =  \int_\Gamma \chi h_e^2 \kappa_B |[\![\mathbf{j}]\!]'|^2 dS && \forall \chi \in \mathbb{V}_0(\Omega), \label{Ohmic_subgrid_heating}
\end{align}
and the proof for energy conservation -- i.e. that the sub-grid scale resistive and Ohmic heating terms cancel in the computation of the rate of change of total energy -- is identical to the continuous case in Proposition 1a. Finally, we note that one may be concerned that a sub-grid scale based Ohmic heating term could lead to spurious small scale noise in the thermal field. In the absence of transport stabilization in the thermal field equation, we did indeed find this to sometimes be the case. However, once transport stabilization was included, this noise vanished. The likely reason behind this is that sub-grid scale resistive damping, and therefore its potentially noise inducing Ohmic heating counterpart, as well as the stabilizing term in the thermal field equation are all triggered by the same local dynamical effect, which is given by transport. In particular, at locations where the sub-grid scale Ohmic heating has a large effect, one would expect the same to be true for the thermal field transport stabilization. Since the transport stabilization term in the thermal field equation is neutral with respect to the internal energy, this process can be seen as locally diffusing the grid-scale Ohmic heating induced internal energy.\\ \\
Next to the above remarks on energy transfer, we find that for the rate of change of total helicity, the sub-grid scale resistive term has an effect analogous to the resolved one. This is to be expected from the filtered ideal MHD equations' point of view, which exhibits reconnection as resolved magnetic field features become unresolved. Following the proof of Proposition 3c for a div-conforming magnetic field $\mathbf{B}$ with the modified electric field $\mathbf{E}_\eta$, we find
\begin{equation}
\frac{d}{dt} \langle \mathbf{A}, \mathbf{B} \rangle =  \dots = -2 \left\langle P_{\mathbb{V}_1(\Omega)}(\mathbf{B}), \mathbf{E}_\eta \right\rangle = \dots - 2\int_\Gamma h_e^2 \kappa_B \left[\!\left[P_{\mathbb{V}_1(\Omega)}(\mathbf{B})\right]\!\right]' \cdot [\![\mathbf{j}]\!]' dS,
\end{equation}
while for the continuous resistive MHD equations, we found in Proposition 1c
\begin{equation}
\frac{d}{dt} \langle \mathbf{A}, \mathbf{B} \rangle =  \dots -2\langle \mathbf{B}, \eta \mathbf{j} \rangle.
\end{equation}
In other words, the sub-grid type resistivity's effect on the total helicity is analogous to the resolved one. In particular, overall the sub-grid type resistivity has the same physical behavior with respect to wave numbers as the non-discretized one (see Remark 2).\\

\textbf{Alternative interior penalty based formulations}\\ \\
In the modified electric field $\mathbf{E}_\eta$ defined by \eqref{E_eta}, we considered an interior penalty term based on a sub-grid resistivity formulation. However, other interior penalty term based formulations may be considered too, which for instance could be aimed at different behaviors with respect to the equations' structure preserving properties. One such attempt could be given by penalizing the transport term
\begin{equation}
\nabla \times (\mathbf{F}/n \times \mathbf{B}), \label{B_transport_term}
\end{equation}
similarly to how the temperature equation's interior penalty term (third term on LHS of \eqref{ideal_T_ec_ip}) penalizes the scalar temperature transport term
\begin{equation}
\mathbf{F} \cdot \nabla T.
\end{equation}
In particular, since the penalty term would then include the momentum $\mathbf{F}$, it would be possible to reestablish energy conservation via a corresponding term in the momentum equation, rather than via a sub-grid scale Ohmic heating term. This could be considered a more natural setup, recalling that from an energetic point of view, the magnetic field transport term relates to the transfer of magnetic energy to and from kinetic energy. Unfortunately, in practice we found interior penalty terms based on \eqref{B_transport_term} or variations thereof to perform badly in terms of stabilization. In future work, it may be of interest to continue to explore this route of magnetic field transport stabilization.\\ \\
Next to the above considerations on energy transfer, one could also attempt to maintain helicity conservation using a suitable interior penalty term.\\ \\
\textbf{Definition 6.} The electric field $\mathbf{E}_{H_M}$ adjusted to contain a helicity-neutral damping effect is given by
\begin{align}
\left\langle \bs{\omega}, \mathbf{E}_{H_M} \right\rangle =& \left\langle \bs{\omega},  - \mathbf{F}/n \times \mathcal{B} \right\rangle + \int_\Gamma h_e^2 \kappa_B \left[\!\left[\mathcal{B} \times \bs{\omega} \right]\!\right]' \cdot \left[\!\left[\mathcal{B} \times \mathbf{j}\right]\!\right]' dS && \forall \bs{\omega} \in \mathbb{V}_i(\Omega), \label{E_AB}
\end{align}
where $k_B$ is a suitable stabilization parameter, and $i=1$ for a curl-conforming electric field, and $i=2$ for a div-conforming one. Further, the magnetic field $\mathcal{B}$ is given by $\mathbf{B}$ when $\mathbf{E}_{H_M}$ is paired with \eqref{BHdiv_discr} or \eqref{BHcurl_discr}, and $P_{\mathbb{V}_1(\Omega)}(\mathbf{B})$ when $\mathbf{E}_{H_M}$ is paired with \eqref{BHdiv_discr_hel}. \begin{center}\vspace{-12mm}\end{center}\hrulefill\\

In terms of the conservation properties for \eqref{E_AB} when paired with \eqref{BHdiv_discr_hel}, we proceed as in \eqref{Bdiv_E_eta1} - \eqref{Bdiv_E_eta2}, and obtain a rate of change of total magnetic energy given by
\begingroup
\addtolength{\jot}{2mm}
\begin{align}
\frac{d}{dt} \int_\Omega \frac{|\mathbf{B}|^2}{2\mu_0} \; dx &= \left\langle \frac{\mathbf{B}}{\mu_0}, \pp{\mathbf{B}}{t} \right\rangle = - \left\langle \frac{\mathbf{B}}{\mu_0}, \nabla \times \mathbf{E}_{H_M} \right\rangle = - \langle \mathbf{j}, \mathbf{E}_{H_M} \rangle +  \int_{\partial \Omega_2} \mathbf{g}_2 \cdot \mathbf{E}_{H_M} \; dS\\
&\hspace{-14mm}= \left\langle \mathbf{j},  \mathbf{F}/n \times P_{\mathbb{V}_1(\Omega)}(\mathbf{B}) \right\rangle - \int_\Gamma h_e^2 \kappa_B \left|\left[\!\left[P_{\mathbb{V}_1(\Omega)}(\mathbf{B}) \times \mathbf{j}\right]\!\right]'\right|^2 dS + \int_{\partial \Omega_2} \mathbf{g}_2 \cdot \mathbf{E}_{H_M} \; dS,
\end{align}
\endgroup
and we now have a non-positive definite stabilizing term based on $|[\![P_{\mathbb{V}_1(\Omega)}(\mathbf{B}) \times \mathbf{j}]\!]'|^2$ rather than $|[\![\mathbf{j}]\!]'|^2$ as obtained from $\mathbf{E}_\eta$. As for $\mathbf{E}_\eta$, the stabilizing interior penalty term dissipates magnetic energy, and to reestablish energy conservation, we modify the temperature equation analogously to \eqref{Ohmic_subgrid_heating}, with $P_{\mathbb{V}_1(\Omega)}(\mathbf{B}) \times \mathbf{j}$ instead of $\mathbf{j}$. Finally, we confirm that \eqref{BHdiv_discr_hel} together with $\mathbf{E}_{H_M}$ conserves helicity. To see this, we follow the proof of Proposition 3c, and obtain
\begingroup
\addtolength{\jot}{2mm}
\begin{align}
\frac{d}{dt} \langle \mathbf{A}, \mathbf{B} \rangle = \dots =& -2 \left\langle \mathbf{B}, \mathbf{E}_{H_M} \right\rangle = -2 \left\langle P_{\mathbb{V}_1(\Omega)}(\mathbf{B}), \mathbf{E}_{H_M} \right\rangle \label{AB_discr_cons1}\\
=& -2 \left\langle P_{\mathbb{V}_1(\Omega)}(\mathbf{B}), -\mathbf{F}/n \times P_{\mathbb{V}_1(\Omega)}(\mathbf{B}) \right\rangle \nonumber\\
&- 2\int_\Gamma h_e^2 \kappa_B \left[\!\left[P_{\mathbb{V}_1(\Omega)}(\mathbf{B}) \times P_{\mathbb{V}_1(\Omega)}(\mathbf{B})\right]\!\right]' \cdot \left[\!\left[P_{\mathbb{V}_1(\Omega)}(\mathbf{B}) \times \mathbf{j}\right]\!\right]'dS= 0, \label{AB_discr_cons2}
\end{align}
\endgroup
where both terms in the last equation cancel since
\begin{equation}
P_{\mathbb{V}_1(\Omega)}(\mathbf{B}) \times P_{\mathbb{V}_1(\Omega)}(\mathbf{B}) = 0. \label{BxB=0}
\end{equation}\\ \\
\textbf{Remark 4.} Before moving on to test the above stabilized electric fields in the numerical results section, we briefly mention another possible total helicity preserving stabilized method. We showed above that the interior penalty term with $\mathcal{B} = P_{\mathbb{V}_1(\Omega)}$ in \eqref{E_AB} is neutral with respect to the rate of change of total helicity due to the vanishing cross product \eqref{BxB=0}. Next to this cross-product based setup, it is also possible to instead consider a scaled dot-product. Such a formulation has been considered from a physics point of view e.g. in \cite{boozer1986ohm}, and the corresponding interior penalty formulation leads to an electric field $\mathbf{E}_{H_{M1}}$ given by
\begin{align}
\left\langle \bs{\omega}, \mathbf{E}_{H_{M1}} \right\rangle =& \left\langle \bs{\omega},  - \mathbf{F}/n \times \mathcal{B} \right\rangle + \int_\Gamma h_e^2 \kappa_B \left[\!\left[\frac{\mathcal{B} \cdot \bs{\omega}}{|\mathcal{B}|^2} \right]\!\right] \left[\!\left[\frac{\mathcal{B} \cdot \mathbf{j}}{|\mathcal{B}|^2}\right]\!\right] dS && \forall \bs{\omega} \in \mathbb{V}_i(\Omega), \label{E_AB1}
\end{align}
where the notation is as in Definition 6, and noting that here we consider jumps of a scalar instead of a vector. Following the above argument \eqref{AB_discr_cons1} - \eqref{AB_discr_cons2} for helicity preservation, we now consider the jump of
\begin{equation}
\frac{P_{\mathbb{V}_1(\Omega)}(\mathbf{B}) \cdot P_{\mathbb{V}_1(\Omega)}(\mathbf{B})}{|P_{\mathbb{V}_1(\Omega)}(\mathbf{B})|^2} = \frac{|P_{\mathbb{V}_1(\Omega)}(\mathbf{B})|^2}{|P_{\mathbb{V}_1(\Omega)}(\mathbf{B})|^2} \equiv 1,
\end{equation}
across facets instead of the cross product \eqref{BxB=0}. Unlike the latter cross product, this time the term is not equal to zero. However, its jump across facets is, and so the overall interior penalty term vanishes again. As before, energy conservation can be reestablished by adding a corresponding Ohmic heating term analogously to \eqref{Ohmic_subgrid_heating}, with $P_{\mathbb{V}_1(\Omega)}(\mathbf{B}) \cdot \mathbf{j}/|P_{\mathbb{V}_1(\Omega)}(\mathbf{B})|^2$ instead of $\mathbf{j}$. For the numerical results section below, we overall found a similar qualitative field behavior when \eqref{E_AB1} was used when compared to $\eqref{E_AB}$. We note that for scenarios including volumes of negligible magnetic field strength, the division by $|\mathcal{B}|^2$ can be numerically challenging; however, for many applications of interest -- such as magnetic confinement fusion -- such scenarios are not of interest.\begin{center}\vspace{-12mm}\end{center}\hrulefill\\

Finally, we note that the stabilized electric fields \eqref{E_eta}, \eqref{E_AB} and \eqref{E_AB1} all penalize jumps of functions of the form $f(\mathbf{j})$. Alternatively, it is also possible to penalize $D f(\mathbf{j})$ instead, where $D$ is a differential operator, and which can be seen as a form of hyper-diffusion.
%
%
%
\section{Numerical results} \label{sec_Numerical_results}
Having reviewed the div- and curl-conforming discretizations, and discussed ways of stabilizing magnetic field transport, we move on to presenting numerical results for the discretizations' structure preserving properties, as well as qualitative and quantitative field development. For simplicity, the space discretizations are coupled with an explicit third order strongly structure preserving Runge Kutta (SSPRK3) time discretization. A description of implicit, energy conserving time discretizations for the temperature and pressure based setups is given in Appendix \eqref{app_impl_time}. In future work, we aim to formulate appropriate nonlinear solver and preconditioning strategies for such implicit schemes, which may for instance be based on efforts for incompressible MHD as described e.g. in \cite{laakmann2022structure, laakmann2022augmented, ma2016robust}, or stabilized MHD discretizations such as \cite{tang2022adaptive}.\\ \\
First, we consider a steady magnetic vortex in a periodic domain, which is based on one described in \cite{balsara2004second}, and which we study with and without a constant background velocity advecting the vortex. In particular, we will use this test case to demonstrate the conservation properties and convergence rates, and to examine the qualitative field development depending on the choice of stabilization method. Second, we further study the stabilization methods' effects on the qualitative field development via a scenario in which a strong horizontal flow vortex twists vertically aligned magnetic field lines.\\ \\
The meshes, finite element discretization, and solver were implemented using the automated finite element toolkit Firedrake \cite{rathgeber2016firedrake}, which relies heavily on PETSc \cite{balay2019petsc}.\\

\begin{figure}[ht]
\begin{center}
\includegraphics[width=1.0\textwidth]{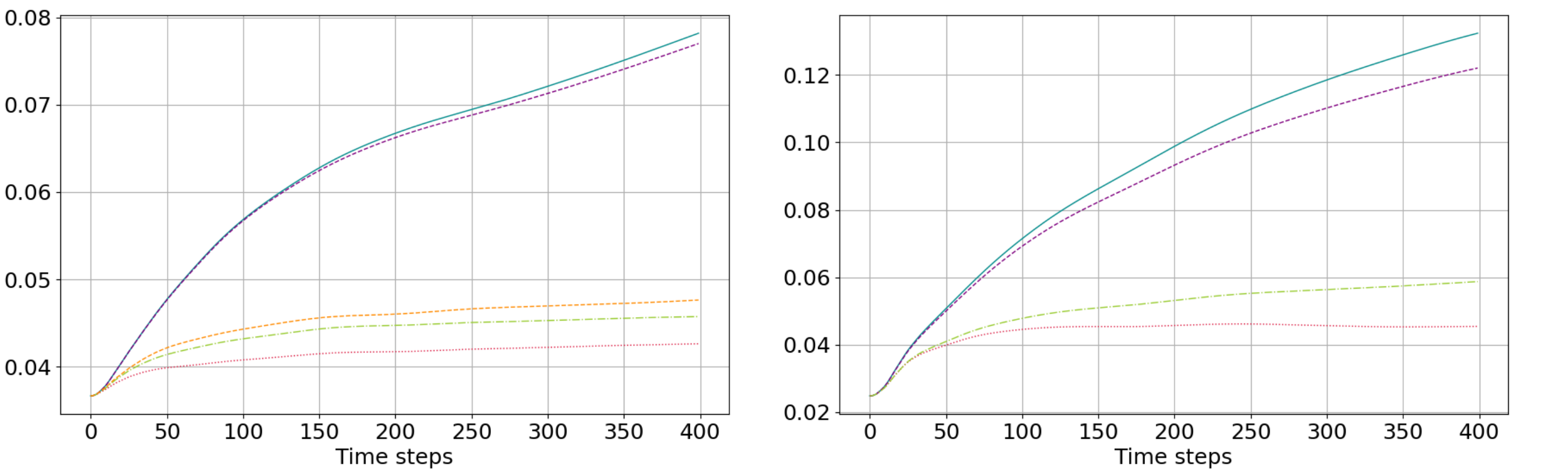}
\caption{Relative error \eqref{L2_error_B} for the magnetic vortex test case with background flow speed $V_b = 0$. Left: div-conforming discretizations. Right: curl-conforming discretizations. The curves correspond to different choices of electric field (see Table \ref{Table_B_discr}), with non-stabilized $\mathbf{E}$ (solid cyan), $\mathbf{E}_{SU\!PG}$ (dashed purple), $\mathbf{E}_\eta$ (dotted red), and $\mathbf{E}_{H_M}$ (dashed green). Dashed orange on the left plot corresponds to the stabilized helicity preserving scheme \eqref{BHdiv_discr_hel} + \eqref{E_AB}.} \label{MV_errors_steady}
\end{center}
\end{figure}

\begin{figure}[ht]
\begin{center}
\includegraphics[width=1.0\textwidth]{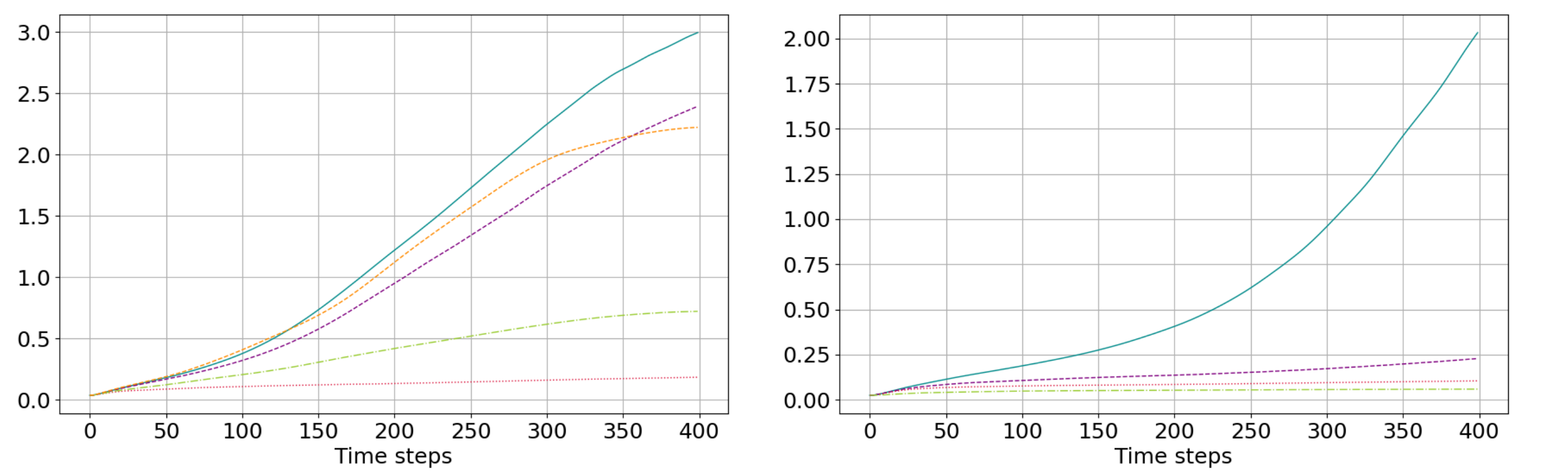}
\caption{Relative error \eqref{L2_error_B} for the magnetic vortex test case with background flow speed $V_b = 0.5$. Arrangement as in Figure \ref{MV_errors_steady}.} \label{MV_errors_nonsteady}
\end{center}
\end{figure}

\begin{figure}[ht]
\begin{center}
\includegraphics[width=1.0\textwidth]{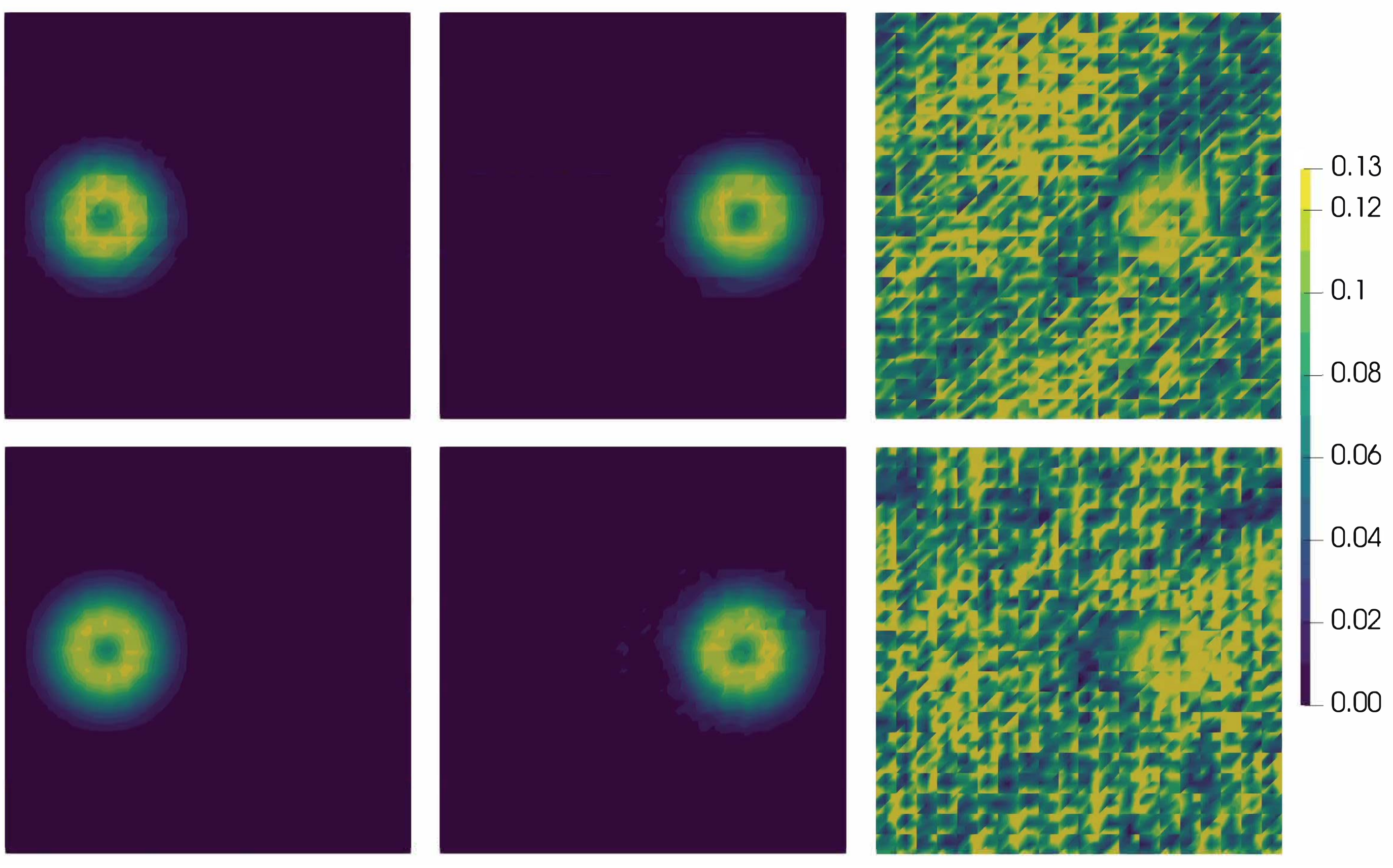}
\caption{Plane $z=2$ for magnetic vortex test case with div-conforming magnetic field (top row) and curl-conforming magnetic field (bottom row), on a periodic regular tetrahedral mesh with finite element chain \eqref{N1_complex}, and background flow speed $V_b = 0.5$. Left column: initial conditions. Center and right columns: simulation at $t=10$, with electric field $\mathbf{E}_\eta$ including sub-grid resistivity-type interior penalty term (center), and without magnetic field transport stabilization (right).} \label{MV_plots}
\end{center}
\end{figure}

\textbf{Magnetic vortex}\\ \\
For the magnetic vortex test case, we consider a periodic domain $\Omega = [0, 10]^2 \times [0, 2]$. The initial conditions are given as functions of
\begingroup
\allowdisplaybreaks
\addtolength{\jot}{2mm}
\begin{align}
&f_{\theta_1}(r) = r e^{1- r^2}, \hspace{1cm} f_{\theta_2}(r) = \int f_{\theta_1}(s)ds \hspace{1cm} f_{\theta_3}(r) = \int \frac{f_{\theta_1}(s)^2}{s} ds\\
&f_{z_1}(r) =
\frac{1}{5r} \begin{cases}
\sin (\pi r) \hspace{36mm} r < \frac{1}{2},\\
\big((1+\kappa)\sin(\pi r) - \kappa \big) \hspace{13mm} r \in \left[\frac{1}{2}, \frac{3}{4}\right],\\
(1 + 2\kappa)(\cos (2 \pi r) - 1)/2 \hspace{6mm} r < 1,\\
0 \hspace{46mm} \text{o.w.},
\end{cases} \hspace{1cm} f_{z_2}(r) \!=\! \int \!s f_{z_1}(s)ds,
\end{align}
\endgroup
where $r = \sqrt{(x - x_c)^2 + (y - y_c)^2}$, denotes the $xy$-plane's radial polar coordinate with respect to the vortex center $(x_c, y_c, 0)^T$, and all integrals are straightforward to evaluate and do not require special functions. Further, $\kappa = (2/\pi - 1/2)/3$, and has been chosen such that $f_{z_2}(r) = 0$ for $r > 1$. For the polar coordinate's angular unit vector $\mathbf{e}_\theta = (-y, x, 0)^T/r$ and vertical unit vector $\hat{\mathbf{z}} = (0, 0, 1)^T$, we then consider initial conditions given by
\begingroup
\allowdisplaybreaks
\addtolength{\jot}{2mm}
\begin{align}
&n = n_0,\\
&\mathbf{V} = V_0 f_{\theta_1} \mathbf{e}_\theta + V_b \hat{\mathbf{x}}, \\
&p = p_0 - \mu_0 \left(f_{\theta_1}^2/2 + f_{\theta_3} + \mu_0 f_{z_1}^2/2\right) + V_0^2 n_0 f_{\theta_3}, \\
&\mathbf{B} = \mu_0 \left(f_{\theta_1} \mathbf{e}_\theta + f_{z_1} \hat{\mathbf{z}}\right),
\end{align}
\endgroup
for $n_0$, $V_0$, $p_0$, $\mu_0$, $m_i$ = $0.5$, $0.5$, $3$, $0.1$, $1$, and background flow speed $V_b$ to be specified below, for a constant background flow in the direction $\hat{\mathbf{x}} = (1, 0, 0)^T$. Finally, the initial magnetic vector potential can be computed as
\begin{equation}
\mathbf{A} = \mu_0(f_{z_2} \mathbf{e}_\phi - f_{\theta_2} \hat{\mathbf{z}}),
\end{equation}
where $\mathbf{A}$ belongs to the same function space as the electric field $\mathbf{E}$. For times $t>0$, we then obtain $\mathbf{A}$ using the same time discretization method as for the full scheme, applied to $\pp{\mathbf{A}}{t} = -\mathbf{E}$. To this end, we note that $\mathbf{A}$ is well defined within the periodic domain since $f_{z_2}$ vanishes for $r>1$. Further, we keep track of $\mathbf{A}$ in order to compute the magnetic helicity $H_M$, which is independent of the scalar potential $\Phi$ appearing in \eqref{A_eqn_cts} due to the periodic domain setup.\\ \\
For the space discretization, we consider a regular tetrahedral mesh and use the discrete de-Rham complex \eqref{N1_complex} corresponding to the second order Nedelec spaces of the first kind. First, we consider the magnetic field's qualitative development depending on the choice of space discretization. For this purpose, we consider two cases given by $x_c$, $y_c$, $V_b$ = $5$, $5$, $0$ and $x_c$, $y_c$, $V_b$ = $2.5$, $5$, $0.5$, and set the space and time resolutions to $\Delta x = 0.5$ and $\Delta t = 0.025$, respectively. Note that we deliberately consider a low space resolution in order to test the magnetic field transport towards under-resolved  scales. The simulations are run up to $t_{max} = 10$, and for each run, we set the temperature field stabilization parameter to $\kappa_T = 0.0001$. Finally, we consider magnetic field discretizations together with discrete electric fields as described in Table \ref{Table_B_discr}. We note that whenever a stabilized electric field is used, we also include the corresponding terms mentioned in the previous section in order to reestablish energy conservation, given by a modified Lorentz force for $\mathbf{E}_{SU\!PG}$, and a sub-grid scale Ohmic heating term for $\mathbf{E}_\eta$ and $\mathbf{E}_{M_H}$.\\
\begin{table}
\begin{center}
\begin{tabular}{|c|c|c|c|} 
 \hline
& $\mathbf{B} \in H(\text{{\normalfont div}};\Omega)$, \eqref{BHdiv_discr} & $\mathbf{B} \in H(\text{{\normalfont div}};\Omega)$, \eqref{BHdiv_discr_hel} &$\mathbf{B} \in H(\text{{\normalfont curl}};\Omega)$, \eqref{BHcurl_discr}\\
 \hline
 $\mathbf{E}$ & $\checkmark$ & $\times$ & $\checkmark$ \\
  $\mathbf{E}_{SU\!PG}$ & $\lambda = 1/4 $& $\times$ &  $\lambda = 1$ \\
  $\mathbf{E}_\eta$ &$\kappa_B = 0.001 $& $\times$ & $\kappa_B = 0.0001$ \\
  $\mathbf{E}_{H_M}$ & $\kappa_B = 0.1$& $\kappa_B = 0.1$ &  $\kappa_B = 0.01$ \\
 \hline
\end{tabular}
\caption{Space discretizations considered in relative magnetic field error analysis for magnetic vortex test case, including choice of magnetic field stabilization parameters.} \label{Table_B_discr}
\end{center}
\end{table}\\
Figure \ref{MV_errors_steady} depicts the relative $L^2$ errors for $\mathbf{B}$ for $V_b = 0$, given by
\begin{equation}
e_{\mathbf{B}}(i) = \frac{\|\mathbf{B}(\mathbf{x}, i\Delta t) - \mathbf{B}_h(\mathbf{x}, i\Delta t)\|_2}{\|\mathbf{B}(\mathbf{x}, i\Delta t)\|_2}, \label{L2_error_B}
\end{equation}
where $i$ denotes the time step counter, and where here the subscript in $h$ denotes the discretized field, while $\mathbf{B}$ denotes the non-discretized expression evaluated at the quadrature points. Further, Figure \ref{MV_errors_nonsteady} depicts the errors for $V_b = 0.5$. Finally, Figure \ref{MV_plots} illustrates the qualitative field development for $V_b = 0.5$, depending on the choice of space discretization\footnote{The magnetic field plots have been generated using Paraview (www.paraview.org). To visualize the div- and curl-conforming finite element fields -- which are based on integral moments -- they are interpolated into a discontinuous Galerkin vector finite element space for tetrahedra; see www.firedrakeproject.org for details. \label{footnote_paraview}}.
\\
\\
We find that without stabilization, there is a clear accumulation of noise in the magnetic field. When the stabilized electric field $\mathbf{E}_\eta$ is used, this noise is suppressed effectively. For $\mathbf{E}_{H_M}$, this is for the most part slightly less so the case, but a clear difference when compared to the non-stabilized runs still exists. On the other hand, the SUPG based electric field discretization $\mathbf{E}_{SU\!PG}$ only leads to a small improvement in all setups except for the curl-conforming discretization when a background velocity is included. For the remaining discussion, we therefore drop the SUPG formulation, and instead only consider the stabilized electric fields $\mathbf{E}_\eta$ and $\mathbf{E}_{M_H}$. \\ \\
\textbf{Remark 4.} One possible explanation for the SUPG stabilization's small error reduction may be that we showed a stabilizing effect in \eqref{SUPG_neg_def} only for the 2D transport term $\tilde{\mathbf{V}} \cdot \nabla B$. However, the original magnetic field transport term, given by
\begin{equation}
\nabla \times (-\mathbf{V} \times \mathbf{B}) = \mathbf{B} (\nabla \cdot \mathbf{V}) - (\mathbf{B} \cdot \nabla)\mathbf{V} + (\mathbf{V} \cdot \nabla)\mathbf{B},
\end{equation}
contains a fully 3D term $(\mathbf{V} \cdot \nabla)\mathbf{B}$, as well as two additional terms that may not be taken into account by the stabilization. Note that the numerical results presented here cannot directly be considered when studying this type of stabilization for its original application area, given by  the vorticity equation. Here, $\mathbf{V}$ and $\mathbf{B}$ are independent prognostic variables, while for the vorticity equation, $\mathbf{V}$ and $\nabla \times \mathbf{V}$ are not.\begin{center}\vspace{-12mm}\end{center}\hrulefill\\

We next demonstrate the space discretizations' energy conservation properties as well as discrete total helicity evolution. For this purpose, we again consider the setups from Table \ref{Table_B_discr}, except for those including $\mathbf{E}_{SU\!PG}$. Further, we include a non-energy conserving setup given by the curl-conforming discretization \eqref{BHcurl_discr} together with $\mathbf{E}_\eta$ \eqref{E_eta}, where, however, the corresponding Ohmic sub-grid heating term \eqref{Ohmic_subgrid_heating} is skipped. The background flow speed is set to $V_b=0$, and we keep the space resolution $\Delta x = 1/20$ used above. The simulations are run up to $t_{max} = 2$, with time steps $\Delta t \in [1/40, 1/80, 1/160, 1/320]$. Results are given in Tables \ref{Table_E_conv} and \ref{Table_MH_conv}, where the values indicate the relative error at the simulations' last time step.
\begin{table}
\begin{center}
\begin{tabular}{|c|c|c|c|c|c|c|c|c|}
 \hline
$\Delta t$ & $\mathbb{V}_2$, $\mathbf{E}$ & $\mathbb{V}_2$, $\mathbf{E}_\eta$ & $\mathbb{V}_2$, $\mathbf{E}_{H_M}$ & $\mathbb{V}_2$, $\mathbf{E}_{H_M}$ $\!\!^\dagger$ & $\mathbb{V}_1$, $\mathbf{E}$ & $\mathbb{V}_1$, $\mathbf{E}_\eta$ & $\mathbb{V}_1$, $\mathbf{E}_\eta$ $\!\!^{\dagger\dagger}$& $\mathbb{V}_1$, $\mathbf{E}_{H_M}$\\
 \hline
$\frac{1}{40}$ & -4.33e-8 & -1.79e-7 & -1.15e-7 & 1.36e-7 & -1.52e-8 & -1.49e-8 & -1.33e-6 & -1.50e-8\\
$\frac{1}{80}$ & -2.82e-8 & -1.37e-7 & -8.18e-8 & -2.03e-9  & -1.09e-8 & -1.08e-8 & -1.32e-6 & -1.08e-8\\
$\frac{1}{160}$ & -8.25e-9 & -6.05e-8 & -3.49e-8 & -1.32e-8 & -3.06e-9 & -3.03e-9 & -1.32e-6 & -3.03e-9\\
$\frac{1}{320}$ & -1.20e-9 & -2.59e-8 & -1.41e-8 & -8.57e-9 & -4.70e-10 & -4.64e-10 & -1.31e-6 & -4.65e-10\\
 \hline
\end{tabular}
\caption{Relative energy errors for steady state magnetic vortex test case, with $t_{max}=10$. $\mathbb{V}_2$ and $\mathbb{V}_1$ indicate div-conforming and curl-conforming magnetic field discretizations \eqref{BHdiv_discr} and \eqref{BHcurl_discr}, respectively. $\mathbf{E}$, $\mathbf{E}_\eta$, and $\mathbf{E}_{H_M}$ indicate electric fields without stabilization, with sub-grid resistive stabilization \eqref{E_eta}, and with helicity preserving sub-grid resistive stabilization \eqref{E_AB}, respectively. Column with $\dagger$ indicates the stabilized fully helicity preserving space discretization $\eqref{BHdiv_discr_hel}$ + \eqref{E_AB}; column with $\dagger\dagger$ indicates run where sub-grid scale Ohmic heating \eqref{Ohmic_subgrid_heating} is skipped.} \label{Table_E_conv} 
\end{center}
\end{table}\\
\begin{table}
\begin{center}
\begin{tabular}{|c|c|c|c|c|c|c|c|c|}
 \hline
$\Delta t$ & $\mathbb{V}_2$, $\mathbf{E}$ & $\mathbb{V}_2$, $\mathbf{E}_\eta$ & $\mathbb{V}_2$, $\mathbf{E}_{H_M}$ & $\mathbb{V}_2$, $\mathbf{E}_{H_M}$ $\!\!^\dagger$ & $\mathbb{V}_1$, $\mathbf{E}$ & $\mathbb{V}_1$, $\mathbf{E}_\eta$ & $\mathbb{V}_1$, $\mathbf{E}_\eta$ $\!\!^{\dagger\dagger}$& $\mathbb{V}_1$, $\mathbf{E}_{H_M}$\\
 \hline
$\frac{1}{40}$ & -1.2e-4 & -1.6e-4 & -1.4e-3 & 1.65e-6 & -1.1e-3 & -1.4e-3 & -1.4e-3 & -1.3e-3 \\
$\frac{1}{80}$ & -1.2e-4 & -1.6e-4 & -1.4e-3 &  8.48e-7 & -1.1e-3 & -1.4e-3 & -1.4e-3 & -1.4e-3 \\
$\frac{1}{160}$ & -1.2e-4 & -1.6e-4 & -1.4e-3 & 4.21e-7 & -1.1e-3 & -1.4e-3 & -1.4e-3 & -1.4e-3 \\
$\frac{1}{320}$ & -1.2e-4 & -1.6e-4 & -1.4e-3 & 2.09e-7 & -1.1e-3 & -1.4e-3 & -1.4e-3 & -1.4e-3 \\
 \hline
\end{tabular}
\caption{Relative magnetic helicity errors for steady state magnetic vortex test case, with $t_{max}=10$. Setup is as in Table \ref{Table_E_conv}.} \label{Table_MH_conv}
\end{center}
\end{table}\\
For the total energy development, we find a reduction as the time step is reduced for all setups except for the one modified to not include a sub-grid heating term, as well as the one corresponding to the stabilized helicity preserving setup \eqref{BHdiv_discr_hel} + \eqref{E_AB}. The reason for the latter setup's apparent lack of reduction in the energy error lies in jumps in the relative total error development in the initial few time steps, combined with a decay following the jumps. This way, the energy error for the coarser time resolutions may happen to be closer to zero at time $t_{max}$ than for the finer ones. The corresponding energy error evolution plot -- which clearly depicts an overall improvement of the energy error development for finer time resolutions -- is given in Appendix \ref{app_E_plots}.\\ \\
For the total helicity development, as expected, we find a lack of convergence for all setups except for the stabilized helicity preserving one given by \eqref{BHdiv_discr_hel} + \eqref{E_AB}. Overall, in terms of helicity preservation, the div-conforming setups perform better than the curl-conforming ones. However, the values for the curl-conforming ones should be taken with caution, since for these setups, the helicity is computed using a div-conforming magnetic vector potential $\mathbf{A}$ computed using its corresponding evolution equation. As discussed in the paragraph before Section \ref{sec_B_eqn_stab}, for curl-conforming $\mathbf{B}$, the div-conforming magnetic vector potential evolution equation is not compatible with the magnetic field evolution equation, thereby possibly introducing an additional error in the computation of the total magnetic helicity. \\ \\
Finally, we test the spatial convergence rates for some of the schemes considered in the above conservation property analysis. In particular, we drop the ones that one would not want consider in practice. From the three schemes based on the stabilized electric field $\mathbf{E}_{H_M}$, we only consider the one based on the fully helicity preserving setup \eqref{BHdiv_discr_hel}, since $\mathbf{E}_{H_M}$ was designed with helicity preservation in mind in the first place. Further, we drop the non-energy conserving setup where sub-grid scale Ohmic heating was skipped. In terms of the test case, we again consider the scenario with $V_b=0$, and $t_{max} = 2$. The time step is fixed to $\Delta t = 1/160$, and we consider spatial resolutions $\Delta x \in [1/20, 1/40, 1/80]$. The convergence rates are computed from the resulting time-averaged relative $L^2$ errors for $\mathbf{B}$, given by
\begin{equation}
\bar{e}_{\mathbf{B}} = \frac{1}{n_t+1} \sum_{i=0}^{n_t} e_{\mathbf{B}}(i),
\end{equation}
where $n_t = t_{max}/\Delta t = 320$ is the number of time steps. The resulting convergence rates are given in Table \ref{Table_conv}. We find that without stabilization, the grid-scale noise may lead to dynamics that degenerate the convergence rate by an increment of more than $0.5$. When stabilization is included, the loss of convergence is less pronounced. This holds true especially for the sub-grid scale resistivity based stabilization \eqref{E_eta}, which leads to a convergence rate in the order of approximately $1.8$. Note that for our choice of finite element complex, we expect an optimal rate of $2$, which may be reduced by an increment of up to $0.5$ when interior penalty terms are included \cite{burman2006edge}. Overall, this again suggests that the schemes including a sub-grid scale resistivity based stabilization term clearly outperform the non-stabilized schemes.
\begin{center}
\begin{table}
\begin{tabular}{|c|c|c|c|c|c|}
 \hline
 & H(div), no stab. & H(div), $\mathbf{E}_\eta$ & H(div), $\mathbf{E}_{H_M}$ & H(curl), no stab. & H(curl), $\mathbf{E}_\eta$\\
$\Delta x$&\eqref{BHdiv_discr} \!+\! \eqref{BHdiv_E_eqn}& \eqref{BHdiv_discr} \!+\! \eqref{E_eta}&\eqref{BHdiv_discr_hel} \!+\! \eqref{E_AB}&\eqref{BHcurl_discr} \!+\! \eqref{BHcurl_E_eqn}& \eqref{BHcurl_discr} \!+\! \eqref{E_eta}\\
 \hline
 $\frac{1}{20} \rightarrow \frac{1}{40}$ & 1.53 & 1.83 & 1.72 & 1.41 & 1.75 \\
 $\frac{1}{40} \rightarrow \frac{1}{80}$ & 1.27 & 1.76 & 1.51 & 1.49 & 1.86 \\
 \hline
\end{tabular}
\caption{Convergence rate for magnetic vortex test case without background velocity.}  \label{Table_conv}
\end{table}
\end{center}

\begin{figure}[ht]
\begin{center}
\includegraphics[width=1.0\textwidth]{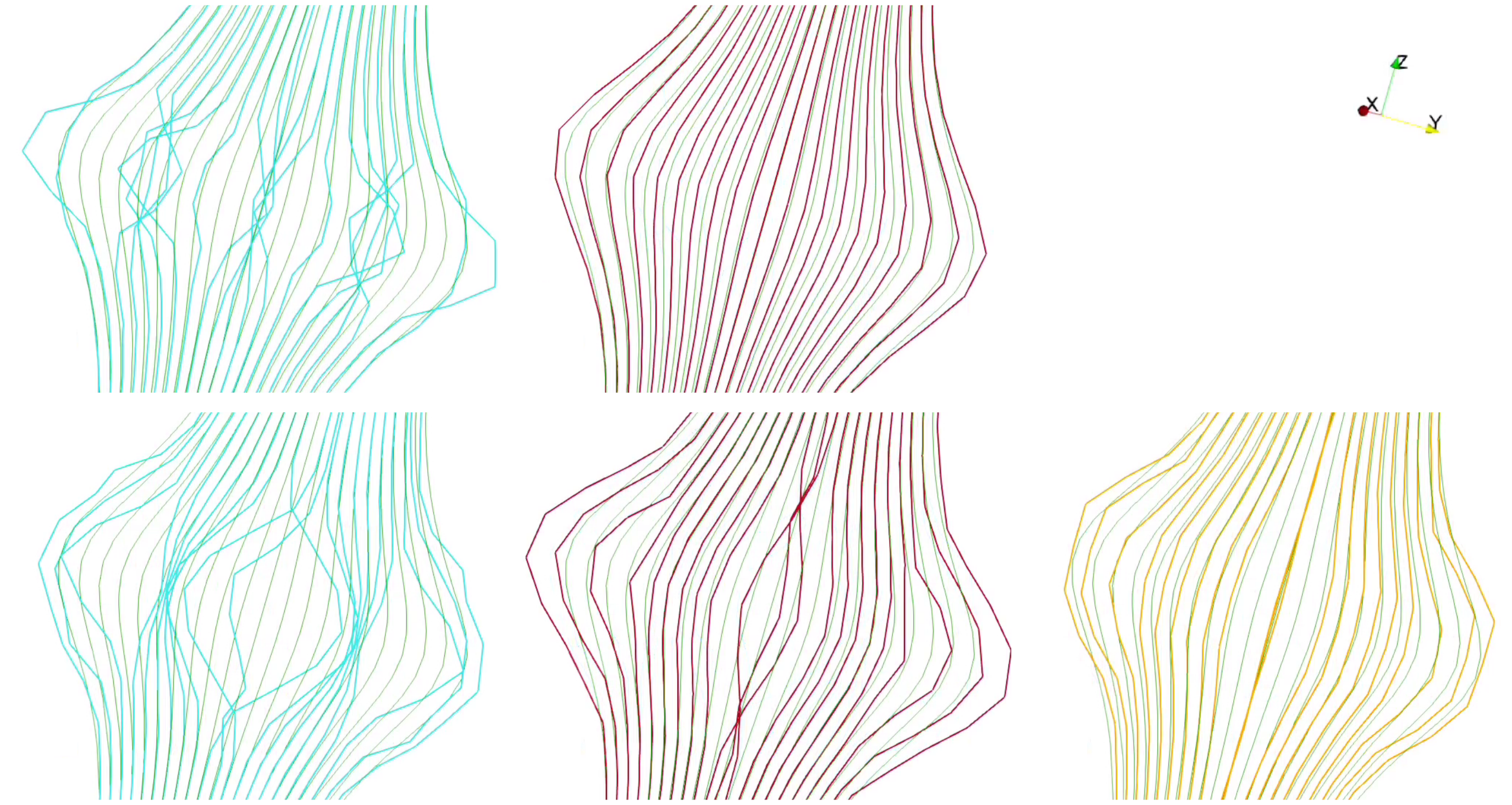}
\caption{Magnetic field lines near vortex center at $t = 20$, for twisted magnetic field test case. Top: curl-conforming discretizations for $\mathbf{B}$. Bottom: div-conforming discretizations for $\mathbf{B}$. The columns correspond to different choices of electric field (see Table \ref{Table_B_discr_twist}), with non-stabilized $\mathbf{E}$ (cyan, left), $\mathbf{E}_\eta$ (red, center), and $\mathbf{E}_{H_M}$ (orange, right). The green field lines in each plot correspond to the reference solution.} \label{twisted_field_line}
\end{center}
\end{figure}

\textbf{Twisted magnetic field lines}\\ \\
Following the above analysis for a steady state test case, we next consider the qualitative magnetic field line development in a scenario consisting of a horizontal flow vortex twisting vertical magnetic field lines. The domain $\Omega$ is given by a cylinder with radius $L_r=5$, and height $L_h = 20$. The initial and boundary conditions are given by
\begingroup
\begin{align}
&n = n_0,\\
&p = p_0,\\
&\mathbf{B} = B_0 \hat{\mathbf{z}}, \hspace{3mm} \mathbf{B} \cdot \mathbf{n}|_{\partial \Omega} = \{-1\; (z=0), \; 0\; (r=L_r), \; 1\; (z=L_h)\},\\
&\mathbf{V} = V_0g(z)\mathbf{e}_\theta
\begin{cases}
0& r< r_1 L_r,\\
\frac{1}{2}\left(1 - \cos\left(\frac{\pi(r - r_1L_r)}{(r_2 - r_1)L_r}\right)\right) & r \in \left[r_1L_r, \; r_2L_r\right],\\
1 &r \in \left(r_2L_r, \; (1-r_2)L_r\right),\\
\frac{1}{2}\left(1 + \cos\left(\frac{\pi(r - (1 - r_2)L_r)}{(r_2 - r_1)L_r}\right)\right) &r \in \left[(1-r_2)L_r, \; (1-r_1)L_r\right],\\
0 & r > (1 - r_1)L_r,
\end{cases}& \mathbf{V} \cdot \mathbf{n}|_{\partial \Omega} = 0,
\end{align}
\endgroup
for $n_0$, $p_0$, $B_0$, $V_0$, = $1$, $2$, $0.1$, $0.5$. Further, we set $g(z) = ((1 - \cos(2\pi z/L_h))/2)^5$, $r_1 = 0.1$, $r_2 = 0.5$, and $\mu_0$, $m_i$ = $3$, $1$. For the space discretization, we consider a regular quadrilateral base mesh (for details, see Appendix \ref{Appendix_mesh}), which is then extruded to create a hexahedral cylindrical mesh. We use the discrete de-Rham complex \eqref{Q1_complex} corresponding to the second order Nedelec spaces of the first kind. The space resolution for the base mesh is given by $\Delta x \approx 0.49$, leading to 180 cells in the $xy$-plane. Further, we divide the vertical extrusion into 32 layers. The time step is given by $\Delta t = 1/32$, and $t_{max} = 20$. For the space discretization, we consider magnetic field setups as described in Table \ref{Table_B_discr_twist}, and the temperature field stabilization is set to $\kappa_T = 0.001$.\\
\begin{table}
\begin{center}
\begin{tabular}{|c|c|c|c|} 
 \hline
& $\mathbf{B} \in H(\text{{\normalfont div}};\Omega)$, \eqref{BHdiv_discr} & $\mathbf{B} \in H(\text{{\normalfont div}};\Omega)$, \eqref{BHdiv_discr_hel} &$\mathbf{B} \in H(\text{{\normalfont curl}};\Omega)$, \eqref{BHcurl_discr}\\
 \hline
 $\mathbf{E}$ & $\checkmark$ & $\times$ & $\checkmark$ \\
  $\mathbf{E}_\eta$ &$\kappa_B = 0.01$& $\times$ & $\kappa_B = 0.001$ \\
  $\mathbf{E}_{H_M}$ & $\times$ & $\kappa_B = 0.1$ &  $\times$ \\
 \hline
\end{tabular}
\caption{Space discretizations considered in qualitative field line analysis for twisted magnetic field test case, including choice of magnetic field stabilization parameters.} \label{Table_B_discr_twist}
\end{center}
\end{table}\\
To assess the qualitative magnetic field development, we examine magnetic field lines\footnote{The field lines are computed in Paraview, using the fields created for visualization purposes (see footnote \ref{footnote_paraview}). They are obtained starting from 20 equidistant points between $(0, -0.5, 0)^T$ and $(0, 0.5, 0)^T$.} at the vortex center's vicinity, near $(r, z) = (0, L_z/2)$. For comparison purposes, we also include a reference run, which is given by the curl-conforming magnetic field discretization including $\mathbf{E}_\eta$, with twice the resolution in space and in time. The corresponding images are given in Figure \ref{twisted_field_line}.\\
\\
We find that as for the magnetic vortex test case, including the interior penalty based magnetic field stabilization methods leads to significant improvements when compared to the non-stabilized setups. For the curl-conforming setup, the electric field $\mathbf{E}_\eta$ given by \eqref{E_eta} leads to a very good field line development when compared to the reference solution. On the other hand, for the div-conforming setup, we find that the helicity preserving discretization with an electric field $\mathbf{E}_{H_M}$ given by \eqref{E_AB} performs slightly better for this test case when compared to $\mathbf{E}_\eta$.
%
\section{Conclusion} \label{sec_conclusion}
In this work, we presented structure preserving space discretizations for the ideal compressible MHD equations, based on the compatible finite element method. We first formulated such discretizations for the equations' hydrodynamic part and included novel energetically neutral interior penalty and upwind methods for the temperature and pressure fields, respectively. We then moved on to the main part of this work, which consists in transport stabilization methods for the magnetic field evolution equation, formulated in a divergence- and a curl-conforming way. For this purpose, we incorporated an SUPG based method that has been considered in the context of vorticity equations before. Further, we devised two new interior penalty based methods, one of which directly corresponds to a sub-grid scale based resistivity, while the other was modified to preserve the discrete helicity. All of these methods crucially still satisfy the discrete zero-divergence property for the magnetic field, since they correspond to modifications within the electric field. Further, for all of these methods, we described possible adjustments to reestablish energy conservation for the full system of equations.\\ \\
In numerical tests for a magnetic vortex scenario, we showed that the SUPG based method underperforms in the setting of magnetic field transport. Further, we showed that the two interior penalty based stabilization methods lead to significant improvements of the qualitative field development. Next to the latter development, we also confirmed the interior penalty based methods' energy conservation properties and helicity conservation property for the helicity-preserving penalty setup. Additionally, we studied the discrete magnetic field equations' spatial orders of convergence, and showed that including stabilization led to improved rates. Finally, we considered a test case in which a flow vortex twists magnetic field lines, and came to the same conclusions as for the magnetic vortex scenario's qualitative field development. Altogether, this suggests that the new interior penalty based formulations are an effective means for magnetic field transport stabilization in the context of exactly zero-divergence maintaining discretizations. Next to their effectiveness, they are also simple to implement in finite element libraries which support building kernels for facet integrals.\\ \\
In future work, we aim to formulate space discretizations for other versions of MHD, such as extended MHD or 2-fluid MHD. Further, we will apply the resulting schemes to relevant physical problems, such as magnetic confinement fusion plasma scenarios. Finally, in view of the latter scenarios, we aim to couple the compatible finite element based space discretizations to more computationally efficient time discretizations and solver strategies.\\ \\
\textbf{Acknowledgments}\\
This work has been supported by the U.S. Department of Energy Office of Fusion Energy Sciences and Office of Advanced Scientific Computing Research under the Tokamak Disruption Simulation (TDS) Scientific Discovery through Advanced Computing (SciDAC) project, as well as the Base Theory Program, both at Los Alamos National Laboratory (LANL) under contract No. 89233218CNA000001. The computations have been performed using resources of the National Energy Research Scientific Computing Center (NERSC), a U.S. Department of Energy Office of Science User Facility operated under Contract No. DE-AC02-05CH11231. We would further like to thank Konstantin Lipnikov and Qi Tang for their advice.
%
%
%
%
\appendix
\section{Fully energy conserving implicit schemes} \label{app_impl_time}
Next to the explicit time discretization employed in the numerical results section, it is also possible to pair the space discretizations introduced in this work with energy conserving implicit time discretizations. For this purpose, we follow the framework established in \cite{cohen2011linear} for Hamiltonian systems, as has been done e.g. in \cite{bauer2018energy, wimmer2020density}. The framework relies on casting the space discretized system of equations in a Poisson system form
\begin{equation}
\frac{dF}{dt} = \{F, H\},
\end{equation}
for any functional $F$ of the prognostic variables, the total energy $H$, and an almost-Poisson bracket $\{,.\}$, which is defined by a antisymmetric bilinear form (for details, see e.g. \cite{salmon1988hamiltonian}). Note that the Poisson system is useful to demonstrate energy conservation, since setting $F=H$ and recalling the bracket's antisymmetry immediately implies conservation of energy. Within the bracket, the arguments $F$ and $H$ appear in form of their variational derivatives with respect to the prognostic variables; indeed the Poisson bracket can be computed using the chain rule for $\frac{dF}{dt}$ together with the prognostic variables' evolution equations, as was done for the total energy $H$ in \eqref{Hamiltonian_full}. The framework in $\cite{cohen2011linear}$ then states that fields in the discretized equations which are originally part of the Poisson bracket are discretized according to the midpoint rule. On the other hand, fields that arise as part of the Hamiltonian derivatives (e.g. \eqref{H_variations} for the temperature based setup), are discretized according to a special time integration rule.\\ \\
For the temperature based setup \eqref{nVT_discr}, we find that the above framework does not apply directly, since the discretized equations cannot easily be cast into a Poisson bracket form. This is because the term including $\frac{1}{2}m_i|\mathbf{V}|^2$ and the pressure gradient term are discretized separately. From a Poisson bracket point of view, on the other hand, they both arise from the Hamiltonian derivative in $n$, and therefore should be discretized equally. However, the framework can still serve as a guide to obtain an energy conserving time discretization (up to the choice of magnetic field boundary conditions), leading to a fully energy conserving scheme of the form
\begingroup
\allowdisplaybreaks
\addtolength{\jot}{2mm}
\begin{subequations} \label{Full_scheme_T}
\begin{align}
&\frac{\Delta n}{\Delta t} + \nabla \cdot \overline{\overline{\mathbf{F}}} = 0, \label{Full_scheme_n_eqn}\\
&\left\langle \mathbf{w}, \frac{\Delta \mathbf{V}}{\Delta t} \right\rangle + \left\langle \nabla \times \left( \mathbf{w} \times \overline{\overline{\mathbf{F}}}/\bar{n} \right), \mathbf{V} \right\rangle - \int_\Gamma \left\{\!\left\{\mathbf{n} \times \left( \mathbf{w} \times \overline{\overline{\mathbf{F}}}/\bar{n} \right)\right\}\!\right\} \cdot \tilde{\mathbf{V}} \; dS \nonumber \\
& \hspace{20mm} - \left\langle \frac{1}{2}\overline{\overline{|\mathbf{V}|^2}}, \nabla \cdot \mathbf{w} \right\rangle + \left\langle \bar{\mathbf{j}}, \frac{\mathbf{w}}{m_i \bar{n}} \times \bar{\mathbf{B}}\right\rangle \nonumber\\ 
&\hspace{20mm} + \left\langle \frac{\mathbf{w}}{m_i \bar{n}}, \nabla ( \bar{n} \bar{T}) \right\rangle - \int_\Gamma [\![ \bar{n} \bar{T}]\!] \left\{\frac{\mathbf{w}}{m_i \bar{n}} \right\}dS = 0 &\forall \mathbf{w} \in \mathring{\mathbb{V}}_2(\Omega), \\
& \left\langle \frac{ \bar{n} \chi}{\gamma -1}, \frac{\Delta T}{ \Delta t}  \right\rangle + \left\langle \frac{ \chi}{\gamma -1} \overline{\overline{\mathbf{F}}}, \nabla \bar{T} \right\rangle + \int_\Gamma \frac{h_e^2 \kappa_T}{(\gamma - 1)\{\bar{n}\}} \left[\!\left[\overline{\overline{\mathbf{F}}} \cdot \nabla \chi\right]\!\right]\left[\!\left[\overline{\overline{\mathbf{F}}} \cdot \nabla \bar{T}\right]\!\right]dS \nonumber \\
& \hspace{2cm} - \left\langle \overline{\overline{\mathbf{F}}}/\bar{n}, \nabla (\chi  \bar{n} \bar{T}) \right\rangle + \int_\Gamma [\![ \chi  \bar{n} \bar{T} ]\!] \left\{ \overline{\overline{\mathbf{F}}}/\bar{n} \right\}dS = 0 & \forall \chi \in \mathbb{V}_0(\Omega),\\
&\left\langle \mathbf{\Sigma}, \frac{\Delta \mathbf{B}}{\Delta t} \right\rangle + \left\langle \nabla \times \mathbf{\Sigma}, - \overline{\overline{\mathbf{F}}}/\bar{n} \times \bar{\mathbf{B}} \right\rangle = 0 & \forall \mathbf{\Sigma} \in \mathbb{V}_1(\Omega),
\end{align}
\end{subequations}
\endgroup
where
\begin{equation}
\frac{\Delta a}{\Delta t} = \frac{a^{n+1} - a^n}{\Delta t}, \hspace{5mm} \bar{a} = \frac{a^{n+1} + a^n}{2},
\end{equation}
for any field $a$, known time level $n$, time level $n+1$ to be computed, and time step $\Delta t$. Finally, the double-barred terms are given by
\begingroup
\addtolength{\jot}{2mm}
\begin{align}
&\overline{\overline{\mathbf{F}}} = P_{\mathring{\mathbb{V}}_2(\Omega)}\left(\frac{1}{3}\left(n^n\mathbf{V}^n + \frac{1}{2}n^n\mathbf{V}^{n+1} + \frac{1}{2} n^{n+1}\mathbf{V}^n + n^{n+1}\mathbf{V}^{n+1}\right)\right),\\
&\overline{\overline{|\mathbf{V}|^2}} = \frac{1}{3}\left(|\mathbf{V}^n |^2 + \mathbf{V}^n \cdot \mathbf{V}^{n+1} + |\mathbf{V}^{n+1} |^2  \right).
\end{align}
\endgroup
Energy conservation can then be demonstrated by computing
\begin{align}
&H(n^n, \mathbf{V}^n, T^n, \mathbf{B}^n) - H(n^{n+1}, \mathbf{V}^{n+1}, T^{n+1}, \mathbf{B}^{n+1}) = \cdots \\
&\hspace{1cm}= \left\langle m_i\frac{1}{2}\overline{\overline{|\mathbf{V}|^2}} + \frac{\bar{T}}{\gamma -1}, \frac{\Delta n}{\Delta t} \right\rangle + \left\langle m_i\overline{\overline{\mathbf{F}}}, \frac{\Delta \mathbf{V}}{\Delta t} \right\rangle + \left\langle \frac{1}{\gamma - 1}, \frac{\Delta T}{\Delta t} \right\rangle + \left\langle \frac{\bar{\mathbf{B}}}{\mu_0}, \frac{\Delta \mathbf{B}}{\Delta t} \right\rangle
\end{align}
for which we rearranged the resulting terms from the total energies at time levels $n+1$ and $n$ to obtain an expression analogously to \eqref{Hamiltonian_full}. Finally, we substitute for the discrete evolution equations \eqref{Full_scheme_T}, that is we take the inner product of \eqref{Full_scheme_n_eqn} with $m_i\frac{1}{2}\overline{\overline{|\mathbf{V}|^2}} + \frac{\bar{T}}{\gamma - 1}$, and set $\mathbf{w} = m_i\overline{\overline{F}}$, $\chi \equiv 1$, and $\mathbf{\Sigma} = \bar{\mathbf{B}}/\mu_0$. Note that this setup works equally for the magnetic field's div-conforming discretizations given in Definition 2 and Proposition 3c, where again we consider a time discretized flux given by $\overline{\overline{\mathbf{F}}}$, and all other magnetic field related terms at the midpoint in time.\\ \\
In terms of the magnetic field transport stabilization methods, it can further be shown that the interior penalty type ones also conserve energy after time discretization if all occurrences of $\mathbf{B}$ and $\mathbf{j}$ are set at the midpoint in time. For the SUPG formulation, the Poisson bracket framework no longer applies in a straightforward manner due to the additional time derivative term in $\mathbf{B}_{res}$ appearing in the modified Lorentz force, and more work is required to still obtain an energy conserving formulation; for details, see \cite{wimmer2020energy}.\\ \\
Finally, the pressure based space discretization \eqref{nVp_discr} together with any of the magnetic field equation space discretizations mentioned in this work can also be combined with an energy conserving time discretization based on the Poisson bracket framework. This can be achieved using a flux-to-velocity recovery operator $\mathbb{U} \colon \mathring{\mathbb{V}}_2(\Omega) \rightarrow \mathring{\mathbb{V}}_2(\Omega)$ as described in \cite{wimmer2020density}, and the resulting fully energy conserving scheme is given by
\begingroup
\allowdisplaybreaks
\addtolength{\jot}{2mm}
\begin{subequations} \label{Full_scheme_p}
\begin{align}
&\left\langle \phi, \frac{\Delta n}{\Delta t} \right\rangle - \left\langle \nabla \phi, \bar{n} \overline{\overline{\mathbb{U}}} \right\rangle + \int_\Gamma \left[\!\left[\phi \overline{\overline{\mathbb{U}}} \right]\!\right] \tilde{\bar{n}} \; dS = 0, \\
&\left\langle m_i \bar{n} \; \mathbf{w}, \frac{\Delta \mathbf{V}}{\Delta t} \right\rangle \!+\! \left\langle \nabla \!\times\! \left( m_i \bar{n} \; \mathbf{w} \times \overline{\overline{\mathbb{U}}} \right), \bar{\mathbf{V}} \right\rangle \!-\!\! \int_\Gamma\! \left\{\!\left\{\mathbf{n} \times \left( m_i \bar{n}\;\mathbf{w} \times \overline{\overline{\mathbb{U}}} \right)\right\}\!\right\} \!\cdot\! \tilde{\bar{\mathbf{V}}} dS \nonumber \\
& \hspace{20mm} + \left\langle \nabla \overline{\overline{P}}, \bar{n}\mathbf{w} \right\rangle - \int_\Gamma \left[\!\left[\overline{\overline{P}} \mathbf{w}\right]\!\right] \tilde{\bar{n}} \; dS  - \langle \bar{p}, \nabla \cdot \mathbf{w} \rangle + \left\langle \bar{\mathbf{j}}, \mathbf{w} \times \bar{\mathbf{B}}\right\rangle = 0 &\hspace{-1cm}\forall \mathbf{w} \in \mathring{\mathbb{V}}_2(\Omega), \\
&\left\langle \frac{\phi}{\gamma - 1}, \frac{\Delta p}{\Delta t}\right\rangle - \left\langle \frac{1}{\gamma - 1} \nabla \phi, \bar{p} \overline{\overline{\mathbb{U}}} \right\rangle + \int_\Gamma \left[\!\left[\phi\overline{\overline{\mathbb{U}}}\right]\!\right] \tilde{p} \; dS  +  \left\langle \bar{p}, \nabla \cdot \overline{\overline{\mathbb{U}}} \right\rangle = 0 & \forall \phi \in \mathbb{V}_3(\Omega),\\
&\left\langle \mathbf{\Sigma}, \frac{\Delta \mathbf{B}}{\Delta t} \right\rangle + \left\langle \nabla \times \mathbf{\Sigma}, - \overline{\overline{\mathbb{U}}} \times \bar{\mathbf{B}} \right\rangle = 0 & \forall \mathbf{\Sigma} \in \mathbb{V}_1(\Omega),
\end{align}
\end{subequations}
\endgroup
where $\overline{\overline{P}} \in \mathbb{V}_3(\Omega)$ and $\overline{\overline{\mathbb{U}}} \in \mathring{\mathbb{V}}_2(\Omega)$ are defined such that
\begin{align}
\overline{\overline{P}} &= P_{\mathbb{V}_3(\Omega)}\left(m_i\frac{1}{2}\overline{\overline{|\mathbf{V}|^2}}\right), \\
\left\langle \bar{n} \mathbf{w}, \overline{\overline{\mathbb{U}}} \right\rangle \!&=\! \left\langle \mathbf{w}, \overline{\overline{\mathbf{F}}} \right\rangle \!=\! \left\langle \mathbf{w}, \frac{1}{3}\!\left(\!n^n\mathbf{V}^n + \frac{1}{2}n^n\mathbf{V}^{n+1} + \frac{1}{2} n^{n+1}\mathbf{V}^n + n^{n+1}\mathbf{V}^{n+1}\!\right)\! \right\rangle & \forall \mathbf{w} \in \mathring{\mathbb{V}}_2(\Omega). \label{bbU_orig}
\end{align}
To show that this scheme conserves energy, we proceed as for the temperature based setup
\begingroup
\addtolength{\jot}{2mm}
\begin{align}
&H(n^n, \mathbf{V}^n, p^n, \mathbf{B}^n) - H(n^{n+1}, \mathbf{V}^{n+1}, p^{n+1}, \mathbf{B}^{n+1}) = \cdots \\
&\hspace{1cm}= \left\langle \overline{\overline{P}}, \frac{\Delta n}{\Delta t} \right\rangle + \left\langle m_i\overline{\overline{\mathbf{F}}}, \frac{\Delta \mathbf{V}}{\Delta t} \right\rangle + \left\langle \frac{1}{\gamma - 1}, \frac{\Delta p}{\Delta t} \right\rangle + \left\langle \frac{\bar{\mathbf{B}}}{\mu_0}, \frac{\Delta \mathbf{B}}{\Delta t} \right\rangle\\
&\hspace{1cm}= \left\langle \overline{\overline{P}}, \frac{\Delta n}{\Delta t} \right\rangle + \left\langle m_i\bar{n}\overline{\overline{\mathbb{U}}}, \frac{\Delta \mathbf{V}}{\Delta t} \right\rangle + \left\langle \frac{1}{\gamma - 1}, \frac{\Delta p}{\Delta t} \right\rangle + \left\langle \frac{\bar{\mathbf{B}}}{\mu_0}, \frac{\Delta \mathbf{B}}{\Delta t} \right\rangle.
\end{align}
\endgroup
To end this section, we note that the above two pressure and temperature based schemes are fully energy conserving only up to the iterative method used to solve the resulting discrete nonlinear system of equations. Typically, only a small number of iterations is performed at each time step in order to save computational cost. In particular, for more nonlinear regimes (e.g. exhibiting turbulence), this often introduces a significant error in energy conservation (see e.g. \cite{lee2021petrov, wimmer2020SUPG}). The latter error may be larger than the error resulting from simplifications to the fully energy conserving scheme, thereby justifying such simplifications if a larger number of nonlinear iterations is computationally too expensive. For instance, we may replace \eqref{bbU_orig} by
\begin{equation}
\left\langle \bar{n} \mathbf{w}, \overline{\overline{\mathbb{U}}} \right\rangle = \left\langle \mathbf{w}, \bar{n}\bar{\mathbf{V}} \right\rangle \hspace{1cm} \forall \mathbf{w} \in \mathring{\mathbb{V}}_2(\Omega) \;\;\; \implies \;\;\; \overline{\overline{\mathbb{U}}} \equiv \bar{\mathbf{V}},
\end{equation}
thereby avoiding the additional computational cost associated with obtaining $\overline{\overline{\mathbb{U}}}$ from \eqref{bbU_orig} at each nonlinear iteration.
\section{Relative energy error plots} \label{app_E_plots}
To clarify the values for the relative energy error associated with the stabilized, helicity preserving magnetic field discretization \eqref{BHdiv_discr_hel} + \eqref{E_AB}, we include the corresponding time evolution plot here. Overall, it can clearly be seen that a reduction in the time step $\Delta t$ leads to a reduction in relative energy error.
\begin{figure}[ht]
\begin{center}
\includegraphics[width=1.0\textwidth]{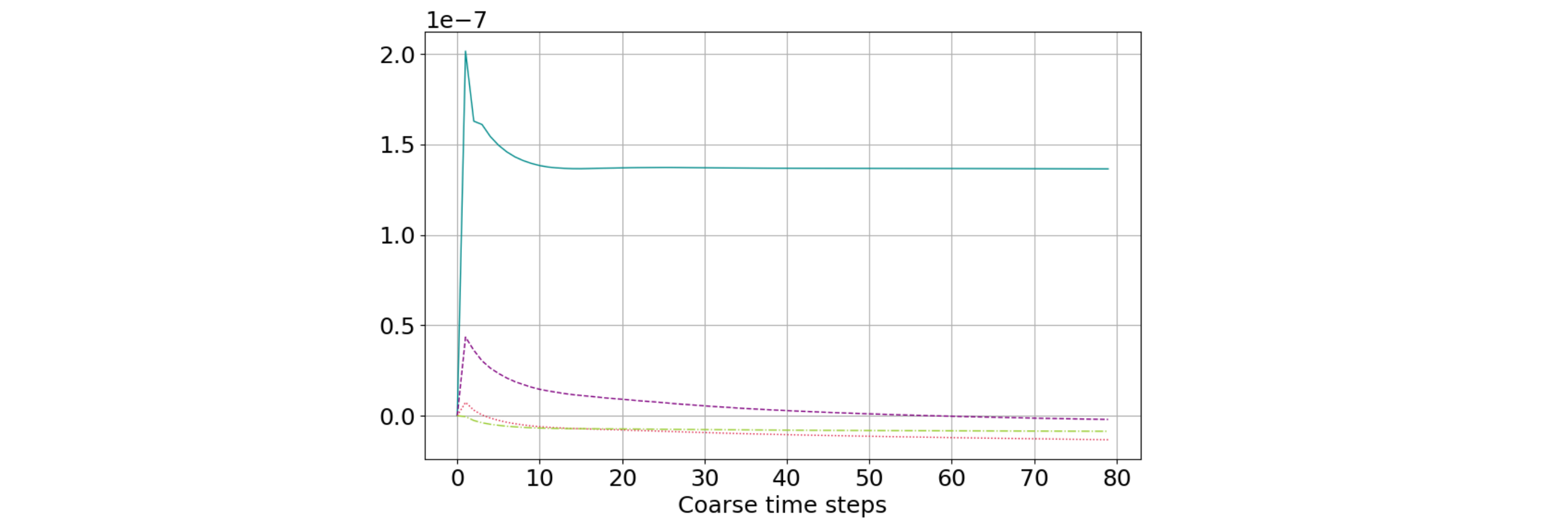}
\caption{Relative energy development for magnetic field space discretization \eqref{BHdiv_discr_hel} + \eqref{E_AB}, for steady magnetic vortex test case. Time resolutions $\Delta t$ given by $1/40$ (solid cyan), $1/80$ (dashed purple), $1/160$ (dotted red), and $1/320$ (dash-dotted green).}
\end{center}
\end{figure}
\newpage
\section{Regular quadrilateral base mesh} \label{Appendix_mesh}
\begin{wrapfigure}{r}{6.6cm}
\begin{center}
\vspace{-6mm}
\includegraphics[width=0.35\textwidth]{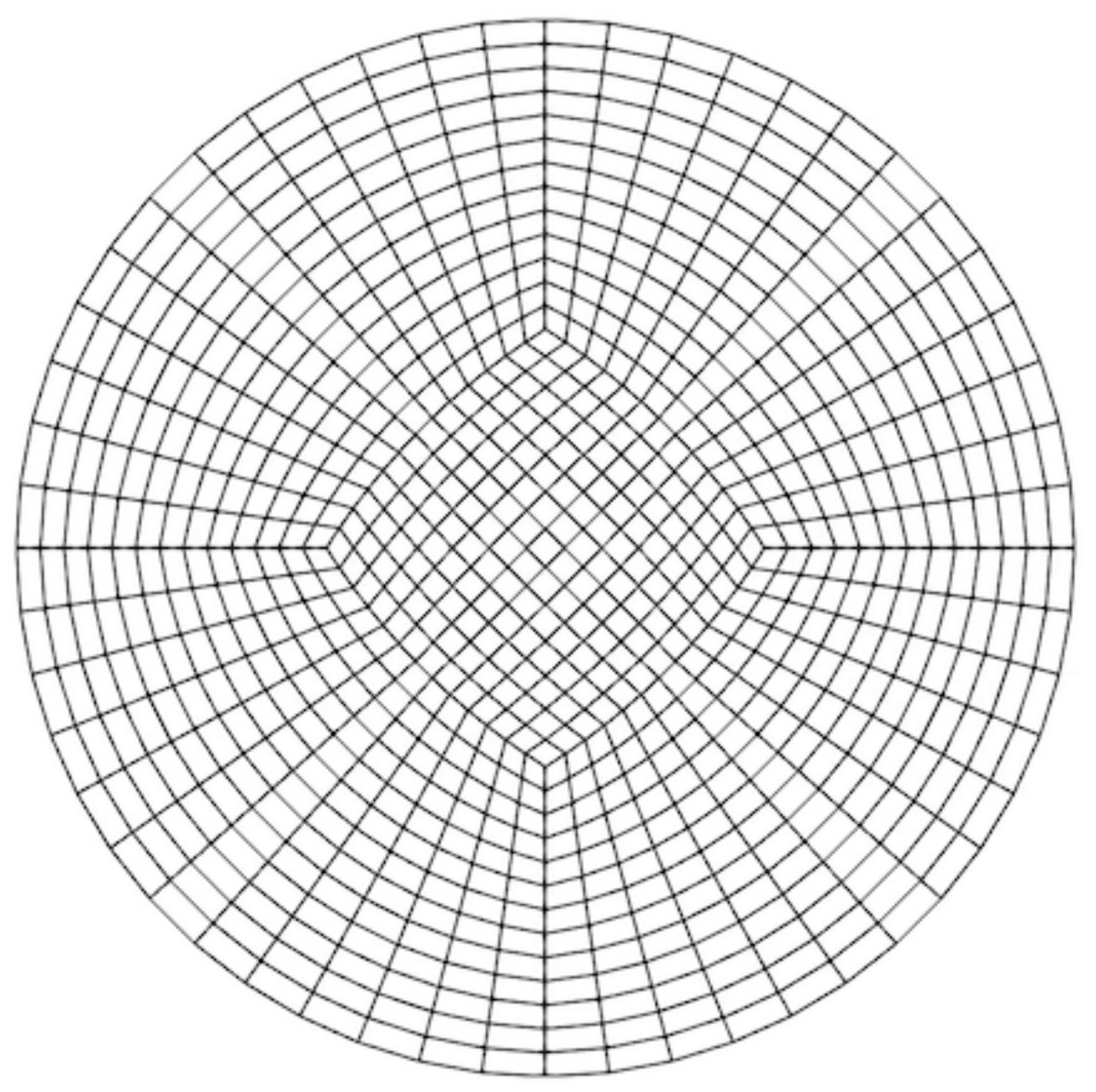}
\end{center}
\vspace{-6mm}
\caption{Base mesh corresponding to $\Delta x = L_r/22$, for radius $L_r$.}
\end{wrapfigure} 
Here, we describe the base mesh used for the extruded mesh in the twisted magnetic field lines test case. It consists of a regular quadrilateral square mesh that is ``bloated'' outward, each side of which is then extruded curvilinearly outward in the $xy$-plane to obtain a disc mesh. The square's side length is equal to $L_s=(2 - \sqrt{2})L_r$, and is chosen such that the latter is equal to the length of the extrusion starting from the square's corners. The square's ``bloated'' feature ensures evenly shaped cells at the square's corners. Overall, given a specified resolution $\Delta x$, we have $5 (L_s/\Delta x)^2$ cells in the base mesh, corresponding to the center square's cells and the ones for the outward extrusions of the four square's edges.

%
%
\bibliographystyle{plain}
\bibliography{Structure_preserving_MHD}
\end{document}